\documentclass[12pt]{amsart}

\usepackage{times}
\usepackage{amsmath,  amssymb,slashed,url,bm,upgreek,amsthm}
\usepackage{graphicx,enumerate}
\newcounter{intro}

\newtheorem{theo}[intro]{Theorem}
\newtheorem{coro}[intro]{Corollary}
\newtheorem{propo}[intro]{Proposition}

\newtheorem{thm}{Theorem}[section]
\newtheorem{lem}[thm]{Lemma}
\newtheorem{prop}[thm]{Proposition}
\newtheorem{cor}[thm]{Corollary}
\newtheorem{defi}[thm]{Definition}
\newtheorem{rem}[thm]{Remark}

\newtheorem*{merci}{Acknowledgements}

\newcommand{\cref}[1]{Corollary~\ref{#1}}

\newcommand{\lref}[1]{Lemma~\ref{#1}}
\newcommand{\pref}[1]{Proposition~\ref{#1}}
\newcommand{\rref}[1]{Remark~\ref{#1}}
\newcommand{\tref}[1]{Theorem~\ref{#1}}

\DeclareMathOperator{\un}{\mathbf{1}}

\DeclareMathOperator{\Id}{Id}

\DeclareMathOperator{\capa}{cap}
\DeclareMathOperator{\card}{card}

\DeclareMathOperator{\ricci}{Ricci}

\DeclareMathOperator{\diam}{diam}
\DeclareMathOperator{\Osc}{Osc}

\DeclareMathOperator{\vol}{dvol}
\DeclareMathOperator{\volu}{vol}
\def\R{\mathbb R}\def\N{\mathbb N}
\def\cB{\mathcal B}
\def\cC{\mathcal C}
\def\cG{\mathcal G}
\def\cL{\mathcal L}
\def\cM{\mathcal M}
\def\cO{\mathcal O}
\def\cS{\mathcal S}
\def\cV{\mathcal V}
\def\cU{\mathcal U}
\begin{document}
\title{Riesz transform on manifolds with quadratic curvature decay}
\author{ Gilles Carron }
\address{Laboratoire de Math\'ematiques Jean Leray (UMR 6629), Universit\'e de Nantes, 
2, rue de la Houssini\`ere, B.P.~92208, 44322 Nantes Cedex~3, France}
\email{Gilles.Carron@univ-nantes.fr}
\begin{abstract}
We investigate the $L^p$-boundedness of the Riesz transform on Riemannian manifolds whose Ricci curvature has  quadratic decay. Two criterions for the $L^p$-unboundedness of the Riesz transform are given. We recover known results about manifolds that are Euclidean or conical at infinity. \\
\\
R\'ESUM\'E:  On \'etudie la continuit\'e de la transform\'ee de Riesz sur les espaces $L^p$  pour des vari\'et\'es dont la courbure de Ricci d\'ecroit quadratiquement. Nous donnons aussi deux crit\`eres g\'eom\'etriques impliquant la non continuit\'e de la transform\'ee de Riesz. Notre m\'ethode nous permet de retrouver les r\'esultats connus pour les vari\'et\'es euclidiennes ou coniques \`a l'infini.
\end{abstract}
\date\today
\maketitle

\footnotetext{\emph{Mots cl\'es: } Transform\'ee de Riesz, fonctions harmoniques, noyau de la chaleur.}

\section{Introduction}
Let $(M^n,g)$ be a complete Riemannian manifold with infinite volume and let $\Delta$ be its associated Laplacian. The Green formula:
$$\forall f\in \cC_0^\infty(M)\,\, ,\, \int_M |df|^2_g \vol_g=\langle \Delta f,f\rangle_{L^2}=\int_M \left|\Delta^{\frac12}f\right|^2 \vol_g\,,$$
implies that the Riesz transform $$R:=d\Delta^{-\frac12}\colon L^2(M)\rightarrow L^2(T^*M)$$ is a bounded operator. It is well known \cite{St2} that on an Euclidean space, the Riesz transform has a bounded extension $R\colon L^p(\R^n)\rightarrow L^p(T^*\R^n)$ for every $p\in (1,+\infty)$. In general,  it is of interest to figure out the range of $p$ for which the Riesz transform extends to a bounded operator $R\colon L^p(M)\rightarrow L^p(T^*M)$ (\cite{St}).
A first remarkable result was obtained by D. Bakry: 
\begin{thm}(\cite{Bakry})\label{thm:Bakry} On a manifold with non negative Ricci curvature, the Riesz transform is bounded on $L^p$, for all $p\in (1,\infty)$. 
\end{thm}
Recall that if $(M,g)$ is a complete Riemannian manifold, its heat kernel $\{h_t(x,y)\}$  is the Schwartz kernel of the heat operator $e^{-t\Delta}$:
$$e^{-t\Delta}f(x)=\int_M h_t(x,y)f(y)dy. $$
According to a well known result of P. Li and S-T. Yau (\cite{LY}), the non negativity of the Ricci curvature implies that  the heat kernel satisfies the upper bound:
\begin{equation}\tag{DUE}\mathrm{for\ all\ } t>0,\ x,y\in M:\ \  h_t(x,y)\le \frac{C}{\volu B(x,\sqrt{t})},\end{equation}
and moreover the Bishop-Gromov inequality implies that the manifold is doubling: there is a constant $\upvartheta$ such that for any $x\in M$ and radius $R>0$:
\begin{equation}\tag{D} \volu B(x,2R)\le  \upvartheta\volu B(x,R)\ .\end{equation}

Another important result is due to T. Coulhon and X-T. Duong: 
\begin{thm}(\cite{CD})\label{thm:CD} The conditions (DUE) and (D)  imply that the Riesz transform is bounded on $L^p$ for any $p\in (1,2]$. \end{thm} 

Manifolds with conical ends satisfy the above conditions (D) and (DUE) and according to the work of H-Q. Li and C. Guillarmou and A. Hassell, we have a complete understanding of the boundedness of the Riesz transform for these Riemannian manifolds. 
\begin{thm}[\cite{Li,GH1, GH2}]\label{Conic} Assume that  $(M^n,g)$ has conical ends: there is a compact set $K\subset M$ such that $(M\setminus K,g)$ is  isometric to a truncated cone:
$$\cC_R(\Sigma):=\,\left((R,\infty)\times\Sigma, (dr)^2+r^2h\right)$$ where 
$(\Sigma,h)$ is a compact Riemannian manifold. 
 
\begin{itemize} \item If $\Sigma$ is not connected, then the Riesz transform is bounded on $L^p$ if and only if $p\in (1,n)\cup\{2\}.$
\item If $\Sigma$ is connected, let $\beta(\beta+n-1)$ with $\beta>0$ be  the first non zero eigenvalue of the Laplacian on $(\Sigma,h)$ and  $\alpha=\min\{\beta,1\}$
then the Riesz transform is bounded on $L^p$ if and only if $p(1-\alpha)<n.$\end{itemize}
\end{thm}

In this paper, we study the boundedness of the Riesz transform on Riemannian manifold whose Ricci curvature satisfies a quadratic decay lower bound, that is to say: 
\begin{equation}\tag{QD} \ricci\ge - \frac{\kappa^2}{r^2(x)} g\,\,,\end{equation}
where  $o\in M$ is a fixed point and $r(x):=d(o,x)$.\par
Manifolds with conical ends satisfy this condition and  our results could be applied to prove the above \tref{Conic} (see subsection \ref{explconic}). Moreover we will also be able to study more general model manifolds (see subsection \ref{WPsection}): 
\begin{thm}\label{warped} Let  $(M^n,g)$ be a  Riemannian manifold  and assume that outside a compact set $(M,g)$ is isometric to the warped product
$$( [1,\infty)\times \Sigma, (dr)^2+r^{2\gamma}h)\, ,$$ where $(\Sigma,h)$ is a compact manifold with non negative Ricci curvature and $\gamma\in (0,1)$.  
\begin{itemize}
\item If $\Sigma$ is connected then the  Riesz transform is bounded on $L^p$ for every $p\in (1,+\infty)$.
\item If $\Sigma$ is not connected then the Riesz transform is bounded on $L^p$ if and only if $1<p\le 2$ or $1<p< (n-1)\gamma+1.$
\end{itemize}
\end{thm}

Our analysis relies on recent results of A. Grigor'yan and L. Saloff-Coste \cite{GS0} and V. Minerbe \cite{Min}  and  we will also use some ideas of P. Auscher, T. Coulhon, X-T. Duong and S. Hofmann \cite{ACDH,AC}.

 We introduce now  two conditions:\begin{itemize}
\item the (VC) (volume comparison) condition: for some constant $C$ and for  any $R\ge 1$ and any $x\in \partial B(o,R)$:
\begin{equation}\tag{VC} \volu B(o,R)\le C\volu B\left(x,R/2\right)\, ,\end{equation}
\item the (RCE) (Relatively Connected to an End) condition: there is a constant $\theta\in (0,1)$ such that for any $R\ge 1$, any $x\in \partial B(o,R)$  there is a continuous path $c\colon [0,1]\rightarrow B(o,R)\setminus B(o,\theta R)$ and a geodesic ray $\gamma\colon [0,+\infty)\rightarrow M\setminus B(o,R)$ satisfying:
\begin{enumerate}
\item $c(0)=x, c(1)=\gamma(0)$, 
\item  the length of $c$ is not too long: 
$\displaystyle \cL(c)\le  \theta^{-1}\, R\,\,,$
\end{enumerate}
The (RCE) condition is an adaption for manifold with several ends of the Relatively Connected Annuli (RCA) condition introduced by  A. Grigor'yan and L. Saloff-Coste in \cite{GS0}.
\end{itemize}
According to \cite{GS0, Min}, under the conditions (QD), (VC) and (RCE), we have a good understanding of the behavior of the heat kernel $\{h_t(x,y)\}$ indeed in this case $(M,g)$ satisfies the doubling condition (D) and the (DUE) estimates (see the discussion in subsection \ref{DDUE}).
Our first result is the following 
\begin{theo}\label{thm:general} Assume that $(M^n,g)$ is a complete Riemannian manifold satisfying the conditions (QD), (VC) and (RCE). If for some positive constants $c$ and $\nu>2$, balls anchored at $o$ satisfy  the reverse doubling hypothesis:
\begin{equation}\tag{RD${}_{\nu}$}\forall R\ge r\ge 1 \,\,\colon \,\,c\left(\frac{r}{R}\right)^\nu\volu B(o,r)\le \volu B(o,R)\, ,\end{equation}   then 
the Riesz transform is bounded on $L^p$ for any $p\in (1,\nu)$.\end{theo}

According to  a beautiful recent result of B. Devyver (\cite[Theorem 5.6]{DevyverP}) our hypothesis on the reverse doubling is equivalent to an isoperimetric inequality for the capacity of anchored balls. Recall that if $\cO\subset M$  is a bounded open subset of a complete Riemannian manifold, then its $p_0$-capacity is defined by:
$$\capa_{p_0} \cO:=\inf\left\{\int_M |d\varphi|^{p_0} \vol_g\,,\, \varphi\in \cC^\infty_0(M) \,\, \mathrm{and}\,\, \varphi\ge 1\,\, \mathrm{on}\,\, \cO\right\}\,.$$ 
And $(M,g)$ is said to be $p_0$-hyperbolic if the $p_0$-capacity of  some/any bounded open subsets  is positive. A non-$p_0$-hyperbolic manifold is called $p_0$-parabolic.  In fact a  Riemannian manifold $(M,g)$ is  $p_0$-parabolic if and only if we can find a sequence of smooth functions with compact support $(\chi_k)$ such that:
$$\left\{\begin{array}{l}
0\le \chi_k\le 1\ ,\\
 \lim_{k\to\infty}\|d\chi_k\|_{L^{p_0}}=0\  \mathrm{and}\\
 \chi_k\to 1, \mathrm{\,\, uniformly\,\, on\,\,  compact \,\, set}.\end{array}\right.$$
 When $p_0=2$, a $2$-hyperbolic manifold is also said to be non-parabolic; moreover a  $2$-hyperbolic manifold  is a manifold carrying a positive Green kernel. A corollary of \tref{thm:general} and of B. Devyver's result is the following: 
 \begin{coro}\label{coro:general} Assume that $(M^n,g)$ is a complete Riemannian manifold satisfying the conditions (QD), (VC) and (RCE). 
If $(M,g)$ is $p_0$-hyperbolic and if the $p_0$-capacity of anchored balls satisfy:
$$\frac{\volu B(o,R)}{R^{p_0}}\le  C\capa_{p_0}(B(o,R))\ ,$$
then the Riesz transform is bounded on $L^p$ for any $p\in (1,p_0)$.\end{coro}

The proof of \tref{thm:general} is based on estimates of the Schwartz kernel of the Riesz transform. Outside the diagonal of $M\times M$, this kernel is smooth and given by:
$$R(x,y)=\int_0^\infty\nabla_xh_t(x,y) \frac{dt}{\sqrt{\pi\, t}}\,.$$
Following the theory of pseudo-differential operators on open manifolds (\cite{Mel1}), which we already used in \cite{CCH}, we separate our analysis in two parts: the part closer to the diagonal  $\{(d(x,y)\le \upkappa r(x)\}$ where we can use the result of \cite{ACDH} and the off diagonal part, where we get the estimate:
$$|R(x,y)|\le \frac{d(x,y)}{r(x)}\frac{C}{\volu B(o,d(x,y))}.$$
When $(M,g)$ is a manifold with Euclidean ends, theses estimates are sharp when one compares to the results of \cite{CCH}.

We will improve an earlier result of  \cite{CCH} and show that if $M$ has two ends then there are very general restrictions on the range of $p$ where the Riesz transform is $L^p$ bounded. \begin{theo}\label{thm:nob} Let $p>2$ and let $(M,g)$ be a $p$-parabolic manifold that is $2$-hyperbolic. If  the Riesz transform is bounded on $L^p$ and on $L^{\frac{p}{p-1}}$ then $M$ has only one end.
\end{theo}
Remark that in this \tref{thm:nob}, no assumptions on the curvature or on the heat kernel is done.

For $q\in (1,2]$, the \tref{thm:CD} provides very general conditions for the $L^q$ boundedness of the Riesz transform  hence, the above criterion is mainly a criterion for the unboundedness of the Riesz transform on $L^p$ when $p>2$. Moreover, this criterion implies that 
the gluing result of B. Devyver is optimal (\cite{Devyver_per}): {\em
 Assume  $M_1,M_2$ are two Riemannian manifolds  that 
satisfy a Sobolev inequality and a lower bound on the Ricci curvature. If on both $M_1$ and $M_2$, the Riesz transform is bounded on $L^p$  and if the connected sum $M_1\#M_2$ is $p$-hyperbolic then
the Riesz transform is $L^p$-bounded on $M_1\#M_2$.}

We recall now several classical notations:\\
{\it \begin{itemize}
\item if $B\subset M$ is a ball we note $r(B)$ its radius and for 
$\theta>0$, $\theta B$ is the ball with the same center and radius $\theta r(B)$.
\item If $\cO\subset M$ and $f\in L^1(\cO)$ we note 
$f_\cO$ the mean of $f$ on $ \cO$:
$$f_\cO=\frac{1}{\volu\cO}\int_{\cO} f.$$
\end{itemize}} 

We recall the following condition about the oscillation of harmonic functions:
\begin{defi} Let $\alpha\in (0,1]$. A complete Riemannian manifold
 $(M^n,g)$ is said to satisfy the scale invariant $\alpha$-H\"older Elliptic (HE${}_{\alpha}$) estimates if there is a constant $C$ such that for any 
ball $B\subset M$ and any harmonic function $h$ defined on $3B$, we have for all $x,y\in B$:
$$|h(x)-h(y)|\le C \,\left(\frac{d(x,y)}{r(B)}\right)^\alpha\, \sup_{z\in 2B} |h(z)|\,\,.$$
\end{defi}
The next result improves \tref{thm:general} in the case where the manifold $M$ has only one end:
 \begin{theo}\label{thm:impro}
Let $(M^n,g)$ be a complete Riemannian manifold with only one end and assume  $(M^n,g)$ satisfies the conditions (QD), (VC) and (RCE) and the reverse doubling hypothesis  (RD${}_{\nu}$) for some exponent $\nu>0$. 

If $(M,g)$ satisfies the scale invariant $\alpha$-H\"older Elliptic estimates, then
the Riesz transform is bounded on $L^p$ for any $p$ such that $p>1$ and $(1-\alpha)p<\nu.$
\end{theo}

Let's explain why this result improves the  \tref{thm:general} when $M$ has only one end.  In the setting of the \tref{thm:impro}, the manifold $(M,g)$ satisfies
the (RCA) condition introduced by A. Grigor'yan and L. Saloff-Coste and we get a scale invariant Poincar\'e inequality: 
for any ball $B\subset M$ and any function 
$f\in \cC^\infty(2B)$ we have
$$\|f-f_B\|_{L^2(B)}\le C r(B)\|df\|_{L^2(2B)}.$$

As $(M,g)$ satisfies the scale invariant Poincar\'e inequalities and the doubling condition (D), the parabolic/elliptic Harnack inequalities hold (\cite{G91,Saloff-Coste}),in particular there is some $\varepsilon\in (0,1]$ such that the property (HE${}_{\varepsilon}$) holds. Hence \tref{thm:impro} yields that the Riesz transform is bounded on $L^p$ as soon as $(1-\varepsilon)p<\nu.$

Our proof is based on a result of P. Auscher and T. Coulhon: assume that $(M,g)$ is a Riemannian manifold satisfying the doubling condition and the scale invariant $L^1$-Poincar\'e inequality. If  there is some $r>p$  such that a $L^{r}$-reverse H\"older inequality holds for the gradient of harmonic functions, then the Riesz transform is bounded on $L^p$ (\cite{AC}, see also \cite{shen}). Although it was not noticed by the authors, this result and the Cheng-Yau's gradient estimate (\cite{ChengYau}) provided another proof of  the \tref{thm:Bakry} of D. Bakry.
 It should be noted (see lemma \ref{Sublne}) that if a manifold $(M,g)$ carries a non constant sublinear harmonic function $h$ with
$$h(x)=\cO\left(r^\beta(x)\right),$$
and if $\alpha>\beta$ then $(M,g)$ can not satisfy the $\alpha$-H\"older Elliptic
 estimates. 
 We will show that the existence of such a sublinear harmonic function yields some restrictions on the range of $p$ where the Riesz transform is $L^p$-bounded:
 \begin{propo}\label{subNO} Let $(M^n,g)$ be a complete Riemannian manifold
 satisfying the conditions (QD), (VC), (RCE) and  the reverse doubling hypothesis (RD${}_{\nu}$) for some exponent $\nu>0$.  Assume moreover that
there are some positive constants $C$ and $\mu$, such that the volume of anchored geodesic balls satisfy: 
 $$\forall R\ge 1\,:\, \volu B(o,R)\le C R^\mu\,\,.$$ If  $(M,g)$ carries a non constant sublinear harmonic function $h$:
$$h(x)=\cO\left(r^\beta(x)\right),$$
then for $p\ge \mu/(1-\beta)$ and $p>\max\left\{\frac{\nu}{\nu-1},2\right\}$,  the Riesz transform can not be bounded on $L^p$.\end{propo}

In  \cite{Cnoharmonic}, we have obtained as a corollary of \tref{thm:impro}:
\begin{coro}\label{noharm}
 Let $(M^n,g)$ be a complete Riemannian manifold that satisfies the curvature quadratic decay  (QD) condition.
If the diameter of geodesic sphere grows slowly
$$\diam \partial B(o,R)= \sup_{x,y\in \partial B(o,R)} d(x,y)=o(R)$$
then the Riesz transform is bounded on $L^p$ for every $p\in (1,+\infty)$.
\end{coro}

In the next section, the condition  (RCE) will be introduced and compared to the (RCA) condition of \cite{GS0}, we will also explain how the analysis of  \cite{GS0, Min} yields the heat kernel estimate (DUE). In the third section, we will explain how the Li and Yau's gradient estimates for solution of the heat equation imply good estimates for the Schwartz kernel of the Riesz transform. The fourth section is devoted to the proof of \tref{thm:general}  and the fifth section to the proof of \tref{thm:impro}. Negative results about the boundedness of the Riesz transform (\tref{thm:nob} and \pref{subNO}) are proved in section 6.
 The theorems \ref{Conic} and \ref{warped} will be proved  in section 7,  we also include examples of manifolds with infinite topological type. Eventually we finish by some perspectives.

\begin{merci} It is a pleasure to thank F. Bernicot, T. Coulhon, B. Devyver, H-J. Hein and E-M. Ouhabaz for useful discussions about my project. I also  thank the referees whose valuable suggestions led to a substantial improvement of the paper. I was partially supported by the ANR grant: {\bf ANR-12-BS01-0004}: {\em Geometry and Topology of Open manifolds}.
\end{merci}
\section{Analysis on manifolds with a quadratic decay of the Ricci curvature}
\subsection{Setting: }In this section, we consider a complete Riemannian manifold $(M^n,g)$ such that for a fixed point $o\in M$, the Ricci curvature satisfies:
\begin{equation}\tag{QD} \ricci\ge - \frac{\kappa^2}{r^2(x)} g\,\,,\end{equation}
where we have defined $r(x)=d(o,x)$. We are going to review geometric conditions that insure that $(M,g)$ satisfies the so-called relative Faber-Krahn inequality: there are positive constants $C,\mu$ such that for any $x\in M$ and $R>0$ and any open domain $\Omega\subset B(x,R)$:
$$\lambda^D_1(\Omega)\ge \, \frac{C}{ R^2}\left(\frac{\volu\Omega}{\volu B(x,R)}\right)^{-\frac2\mu}\,\, .$$
We have noted $\lambda^D_1(\Omega)$ the lowest eigenvalue of the Dirichlet Laplacian on $\Omega$:
$$\lambda^D_1(\Omega)=\inf_{\varphi\in \cC^\infty_0(\Omega)}\frac{\int_\Omega |d\varphi|^2}{\int_\Omega |\varphi|^2}\,\,.$$
Our discussion below is based on the results of A. Grigor'yan and L. Saloff-Coste \cite{GS0} and V. Minerbe \cite{Min}.

It is well known \cite{GriRev} that the relative Faber-Krahn inequality are equivalent to the conjunction of the following two properties:
\begin{itemize}
\item $(M,g)$ is doubling: there is a constant $\upvartheta$ such that for any $x\in M$ and radius $R>0$:
$$\volu B(x,2R)\le  \upvartheta\volu B(x,R).$$
\item The heat kernel $h_t(x,y)$ satisfies the upper bound: for all $t>0$ and all $x,y\in M$ ,
$$h_t(x,y)\le \frac{C}{\volu B(x,\sqrt{t})}.$$
\end{itemize}
According to \cite[Proposition 5.2]{GriRev}, the above relative Faber-Krahn inequality implies (for the same exponent $\mu$) that for some constant $C>0$:
$$\mathrm{for\ all\ }x\in M,\ 0<r<R\ :\ \ \volu B(x,R)\le C\left(\frac{R}{r}\right)^\mu\volu B(x,r) .$$
\subsection{Remote balls} A ball $B(x,\rho)\subset M$ is called remote if its center $x$  and radius $\rho$ satisfy
$$\rho\le \frac{r(x)}{2}.$$ Note that on the remote ball $B(x,r(x)/2)$ the Ricci curvature satisfies:
$$\ricci_g\ge -\frac{4 \kappa^2}{r^2(x)}\, g_x\,\,,$$ hence
the Bishop-Gromov comparison theorem implies that all remote balls satisfy the doubling condition:
If $B$ is a remote ball and if $\theta\in (0,1)$ then
\begin{equation}\label{VDR}
\theta^n \volu(B)\le C(n,\kappa) \volu(\theta B)\,.\end{equation} 
Similarly, the condition (QD) implies that  remote balls satisfy the Poincar\'e inequality (\cite[inequality (4.5)]{buser}) and a  relative Faber-Krahn inequality \cite[Theorem 2]{G91} or \cite[Theorem 3.1]{Saloff-CosteJDG}:
\begin{lem}
If $B\subset M$ is a remote ball, then 
for all $\varphi\in \cC^1(B)$:
$$\left\|\varphi-\varphi_B\right\|^2_{L^1(B)}\le B(n,\kappa)\, r(B)\left\|d\varphi\right\|_{L^1(B)}\,\,,$$ and for all domain $\Omega\subset B$
$$\lambda^D_1(\Omega)\ge \frac{C(n,\kappa) }{ r(B)^2}\, \left(\frac{\volu\Omega}{\volu B}\right)^{-\frac2n}\,\,.$$
\end{lem} 

A ball centered at $o$ will be called {\it anchored}.
\subsection{The doubling condition}According to \cite[Proposition 4.7]{GS0}, this condition  is insured by the doubling of remote balls (\ref{VDR}) and by  the volume comparison (VC) assumption:
there is a constant $C$ such that for any $x\in M$:
\begin{equation}\label{eq:VC}
\tag{VC} \volu B(o,r(x))\le C \volu B\left(x,r(x)/2\right)
\end{equation}

We recall that the doubling condition implies that the volume of balls varies slowly with the center of balls:
for any $\gamma\ge 1$ there is a constant $C_\gamma$ such that if $d(x,y)\le \gamma R$ and $\gamma^{-1}R\le r\le \gamma R $ then\,:
$$C_\gamma^{-1}\le \frac{\volu B(x,R)}{\volu B(y,r)}\le C_\gamma\,\,.$$
In particular if $R\ge r(x)/\gamma$ then 
\begin{equation}\label{eq:comparison}
C_\gamma^{-1} \volu B(o,R)\le \volu B(x,R)\le C_\gamma \volu B(o,R).\end{equation}

Finally it is well known (see \cite[Lemma 2.10]{GS0}) that for a connected, non-compact manifold, the doubling condition implies a reverse doubling condition: there is a constant $\delta>0$, depending only on the doubling constant $\upvartheta$ such that for all $x\in M$ and all $R>r$
$$\delta \left(\frac Rr\right)^\delta\le \frac{\volu B(x,R)}{\volu B(x,r)}.$$
\subsection{Geometry of annuli}

\subsubsection{Number of ends} We remark that the doubling condition implies that $(M,g)$ has a finite of ends: i.e. there is an integer $N$ such that for any $R$, $M\setminus B(o,R)$ has at most $N$ unbounded connected components.
Indeed if $\cO\subset M\setminus B(o,R)$ is an unbounded connected component, there is a point $x_\cO\in \cO\cap\partial B(o,2R)$,
 we have the inclusions $B(x_\cO,R)\subset \cO$ and $B(x_\cO,R)\subset B(o,3R)$, hence we get 
 $$\sum_{\cO} \volu B(x_\cO,R)\le \volu B(o,3R)\, .$$
 However using the doubling condition we get:
 $$\volu B(o,3R)\le \volu B(x_\cO, 5R)\le \upvartheta^3 \volu B(x_\cO,R)\,\,.$$
Hence $M\setminus B(o,R)$ has at most $\upvartheta^{3}$ unbounded connected components.

A slight variation of this argument shows that for any $\lambda>1$, the annulus
$A_{\lambda ,R}=B(o,\lambda R)\setminus B(o,R)$ has at most $N(\lambda,\upvartheta)$ connected components that intersects $\partial B(o,\lambda R)$, however these connected components do not necessary intersect an unbounded connected component of $M\setminus B(o,R)$.
\subsubsection{The (RCE) condition}
In \cite{GS0}, A. Grigor'yan and L. Saloff-Coste have introduced the Relatively Connected Annuli (RCA) condition:
\begin{defi} A manifold $(M,g)$ is said to satisfy to (RCA) condition if there is a point $o\in M$ and a constant $\theta\in (0,1)$ such that for any $R>0$ and any points $x,y\in \partial B(o,R)$ there is a $\cC^1$ path $c\colon [0,1]\rightarrow M$ satisfying:\begin{itemize}
\item  $c(0)=x$, $c(1)=y$,
\item $\displaystyle \cL(c)\le r(x)/\theta,$ 
\item $\displaystyle c([0,1])\subset B(o, \theta^{-1}R)\setminus B(o, \theta R).$ 
\end{itemize}
\end{defi}

It is easy to show that the (RCA) condition implies that $M$ has only one end: i.e. for any compact set $K\subset M$, $M\setminus K$ has only one unbounded connected component. The (RCE) condition is an adaptation of the (RCA) condition for manifolds with several ends.
  
\begin{defi}\label{RCE} We say that a complete Riemannian manifold $(M,g)$ with a finite number of ends satisfies the Relatively Connected to an End (RCE) condition if there is a constant $\theta\in (0,1)$ such that for any point $x$ with $r(x)\ge 1$ there is a continuous path $c\colon [0,1]\rightarrow M$ satisfying:\begin{itemize}
\item  $c(0)=x$,
\item the length of $c$ is bounded by $\cL(c)\le r(x)/\theta,$ 
\item $\displaystyle c([0,1])\subset B(o, \theta^{-1}r(x))\setminus B(o, \theta r(x)),$ 
\item There is a geodesic ray $\gamma:[1,+\infty)\rightarrow M\setminus B(o,  r(x))$ with $\gamma(0)=c(1)$
\end{itemize}
\end{defi}

The (RCE) condition implies that any point can be connected to an end by a path that stays at bounded distance away from the origin. It is easy to see that if $M$ has only one end, the (RCE) condition is just the (RCA) condition of A. Grigor'yan and L. Saloff-Coste.

If $M$ has a finite number of ends, then there is a finite number of geodesic rays $c_1,\dots,c_r\colon [0,+\infty)\rightarrow M$ with $c_i(0)=o$ such that for every $R>>1 $, $M\setminus B(o,R)$ has exactly $r$-unbounded connected component $\cO_1,\dots,\cO_r$ and for all $i$: $$c_i\left((R,+\infty)\right)\subset \cO_i\ .$$
In this setting, we could  replace the last condition in the definition \ref{RCE} by:
\begin{itemize}
\item there is some $i\in\{1,\dots,r\}$ such that $c(1)=c_i(r(x))$.
\end{itemize}

\subsection{Relative Faber-Krahn inequality}\label{DDUE}
The results of A. Grigor'yan and L. Saloff-Coste and of V. Minerbe imply:
\begin{thm}\label{FK}
Assume that $(M^n,g)$ is a complete Riemannian manifold  satisfying the conditions (QD), (VC) and (RCE). Then $(M,g)$ satisfies the relative Faber-Krahn inequality: for some $\mu>0$ and $C>0$ and for any ball $B\subset M$ and any domain $\Omega\subset B$
$$\lambda^D_1(\Omega)\ge \frac{C }{ r(B)^2}\, \left(\frac{\volu\Omega}{\volu B}\right)^{-\frac2\mu}.$$
When $M$ has only one end, then $(M,g)$ satisfies the scale invariant $L^1$ Poincar\'e inequality:
for any ball $B$ and any function $f\in \cC^\infty(2B)$ then 
$$\|f-f_B\|_{L^1(B)}\le C r(B) \, \|df\|_{L^1(2B)}\,\,.$$
\end{thm}

The second assertion is one of the  main result of \cite[thm 5.2]{GS0}--a priori the article deals with the scale invariant $L^2$ Poincar\'e inequality but the argument carries over the case of any $L^p$ Poincar\'e inequality. Strito senso, the first assertion can not  be found in the paper of V. Minerbe; however a quick glimpse on the argumentation shows that the limitation on the exponent $\nu>1$ in the reverse doubling condition,
\begin{equation}\tag{RD${}_{\nu}$}\forall R>r:\,\, \volu B(o,R)\ge\varepsilon  \left(\frac Rr \right)^\nu\volu B(o,r)\end{equation} is made only to insure the (RCA) condition. Under the assumptions of the \tref{FK}, the proof of \cite[Theorem 2.19]{Min} implies that for any $p\ge n$ with $p>2$ there is a constant $C$ such that the weighted Sobolev inequality holds
$$\forall f\in \cC^\infty_0(M):\, \left(\int_M\! |f(x)|^{\frac{2p}{p-2}}dx\right)^{1-\frac2p} \le  \int_M{\displaystyle \frac{C r(x)^2}{(\volu B(o,r(x)))^{\frac2p}} }|df|^2(x)dx.$$

Let's explained how this inequality implies the  relative Faber-Krahn inequality 
for anchored balls i.e. balls centered at $o$:
the doubling condition yields a constant $\mu$ such that 
 $$\forall R>r\,\,:\,\, \volu B(o,R)\le  C \left(\frac Rr \right)^\mu\volu B(o,r)\,\,.$$
 Remark that if this inequality is true for some $\mu$, then it holds for any $\mu'\ge \mu$. Moreover looking at the limit $r\to 0$, we see that $\mu\ge n$ and in the following we will assume that $\mu>2$.
 
In particular using the above Sobolev inequality for $p=\mu$, we get that for any function $f\in \cC^\infty_0(B(o,R))$:
\begin{equation*}\begin{split}
\left(\int_{B(o,R)} |f(x)|^{\frac{2\mu}{\mu-2}}dx\right)^{1-\frac2\mu} 
&\le C \int_{B(o,R)} \frac{r(x)^2}{(\volu B(o,r(x)))^{\frac2\mu}} |df|^2(x)dx\\
&\le C \frac{R^2}{(\volu B(o, R))^{\frac2\mu}}\int_{B(o,R)} |df|^2(x)dx\,\,.
\end{split}\end{equation*}
Then with the H\"older inequality, we get for any domain $\Omega\subset B(o,R)$:
$$1\le C\frac{R^2}{(\volu B(o, R))^{\frac2\mu}} (\volu\Omega))^{\frac2\mu} \lambda^D_1(\Omega)\,\,.$$

It is now easy to show that the  relative Faber-Krahn inequality holds for all balls. Indeed it remains to show the relative Faber-Krahn inequality  for balls $B(x,\rho)$ with $\rho\ge r(x)/2$. But such a ball satisfies:
$$B(x,\rho)\subset B(o,\rho+r(x))\subset B(o,3\rho)\ \mathrm{and}\ B(o,3\rho)\subset B(x,5\rho)\, .$$
Hence if $\Omega\subset B(x,\rho)$ then
\begin{equation*}\begin{split}
\lambda_1^D(\Omega)&\ge \frac{C}{(3\rho)^2} \left(\frac{\volu B(o,3\rho)}{\volu \Omega}\right)^{\frac2\mu}\\
&\ge \frac{C}{(3\rho)^2}\left(\frac{\upvartheta^{-2}\volu B(o,5\rho)}{\volu \Omega}\right)^{\frac2\mu}\\
&\ge  \frac{C}{(3\rho)^2}\left(\frac{\upvartheta^{-2}\volu B(x,3\rho)}{\volu \Omega}\right)^{\frac2\mu}\\
&\ge  \frac{C}{(3\rho)^2}\left(\frac{\upvartheta^{-4}\volu B(x,\rho)}{\volu \Omega}\right)^{\frac2\mu}\ .
\end{split}\end{equation*}
In fact a remarkable result of V. Minerbe \cite[Prop. 2.8]{Min} (see also \cite[Proposition 4.5]{HaK} for an earlier result) shows that the (RCA) condition is insured by an anchored Poincar\'e inequality and a reverse doubling condition:
\begin{thm} Assume that $(M,g)$ is a complete Riemannian manifold that is doubling and such that balls $B=B(o,R)$ centered at $o$ satisfy the Poincar\'e inequalities:
$$\forall f\in \cC^\infty(2B) \,\,;\, \|f-f_{B}\|_{L^p(B)}\le C R \|df\|_{L^p(2B)}\ .$$ If for positive constants $C$ and $\nu>p$,  we have the reverse doubling condition: 
$$\forall R>r\,\,:\,\, \volu B(o,R)\ge C \left(\frac Rr \right)^\nu\volu B(o,r)$$
then $(M,g)$ satisfies the (RCA) condition.
\end{thm}

\section{Estimates on the Riesz kernel}
\subsection{Assumptions}
In this section, we assume that $(M^n,g)$ is a complete Riemannian manifold with a based point $o\in M$ satisfying the following conditions:
\begin{enumerate}[i)]
\item A quadratic decay  on the negative part of the Ricci curvature 
$$\ricci_g\ge -\frac{\kappa^2}{r^2(x)}\, g.$$
\item There are positive constants  $\mu$ and $C$, such that  the relative Faber-Krahn inequality holds: for any ball $B\subset M$ and any domain $\Omega\subset B$
$$\lambda^D_1(\Omega)\ge \frac{C }{ r(B)^2}\, \left(\frac{\volu\Omega}{\volu B}\right)^{-\frac2\mu}.$$
\item For some positive constants $c$ and $\nu>2$, we have the reverse doubling condition (RD${}_{\nu}$) for anchored balls:
$$\forall R>r\,\,:\,\, c  \left(\frac Rr \right)^\nu\volu B(o,r)\le \volu B(o,R) .$$
\end{enumerate}
\begin{rem} The limitation $\nu>2$ is not essential, we can handle the case where $\nu>1$ as well but in the case $\nu\in (1,2]$, the estimate on the Riesz kernel is more complicated and the conclusion of the main theorem are interesting only when $\nu>2$. Indeed the relative Faber-Krahn hypothesis implies that if  $p\in (1,2]$, then the Riesz transform is bounded on all $L^p$ (\tref{thm:CD}   or \cite{CD}).
\end{rem}
\subsection{Li and Yau's inequality}When $B\subset M$ is a remote ball, then on $\frac32 B$ the Ricci curvature is bounded from below by $-16\kappa^2 r(B)^{-2}$, so that according to P. Li and S-T. Yau's Harnack inequality \cite[Theorem 2.1]{LY}: there is a constant $c(n,\kappa)$ such that for any positive solution of the heat equation $u\colon [0,2T]\times \frac32 B\rightarrow \R_+^*$, we have on $[T,2T]\times B$:
\begin{equation}\label{LYG}
\frac{|\nabla u|^2}{u^2}-2\frac1u \frac{\partial u}{\partial t}\le c(n,\kappa) \left(\frac1T+\frac{1}{r^2(B)} \right).\end{equation}
\subsection{Spatial derivative of the heat kernel}
According to A. Grigor'yan \cite{GriRev,G1}, in our setting,the heat kernel  satisfies the following Gaussian upper bound: for all $t>0$ and $x,y\in M$, we have 
$$h_t(x,y)\le \frac{C}{\volu B(x,\sqrt{t})} e^{-\frac{d^2(x,y)}{ct}}\,;$$
$$\left|\frac{\partial}{\partial t}h_t(x,y)\right|\le \frac{C}{\sqrt{t} \volu B(x,\sqrt{t})} e^{-\frac{d^2(x,y)}{ct}}.$$
Let $t>0$ and $x,y\in M$; on the parabolic ball $ [0,t/2]\times B(x,\sqrt{t})$, the function $u(s,z):=h_{\frac t2+ s}(z,y)$
 satisfies:
\begin{equation}\label{GaussG2}u(s,z)+\sqrt{t} \left|\frac{\partial}{\partial s} u(s,z)\right|\le \frac{C}{ \volu B(x,\sqrt{t})} e^{-\frac{d^2(x,y)}{ct}}.\end{equation}
Introduce now $\rho=\min\left\{\sqrt{t}, r(x)/2\right\}$, the ball $B(o,\rho)$ is remote and  Li and Yau's above estimate (\ref{LYG}) yields the following:
$$ |\nabla u|^2(t/2,x)\le 2u(t/2,x)\frac{\partial u}{\partial t}(t/2,x)+C \left(\frac{1}{t}+\frac{1}{\rho^2} \right)u^2(t/2,x)$$ 
and with the estimate (\ref{GaussG2}), we get: \begin{equation*}
\left|\nabla_x h_t(x,y)\right|\le \left(\frac{1}{\sqrt{t}}+\frac{1}{r(x)}\right) \frac{C}{ \volu B(x,\sqrt{t})} e^{-\frac{d^2(x,y)}{ct}}\,.
\end{equation*}

\subsection{Application to the Schwartz kernel of the Riesz transform}
Recall that the Riesz transform is the operator
$$R=d\Delta^{-\frac12}\colon L^2(M)\rightarrow L^2(T^*M),$$
its Schwartz kernel is smooth on $M\times M\setminus\mathrm{Diag}$:  if $x\not=y\in M$ then $R(x,y)\in T_x^*M$  is given by:
$$R(x,y)=\int_0^{+\infty}\nabla_xh_t(x,y)\frac{dt}{\sqrt{\pi \,t\,}}\,\, .$$

Let $\upkappa\ge 4$, we are going to estimate $|R(x,y)|$ in three different regimes:
\begin{enumerate}[i)]
\item First regime: $d(x,y)\ge \frac{1}{\upkappa} r(x)$ and $\frac{1}{\upkappa} r(x)\le r(y)\le \upkappa r(x)$,
\item The short to long range regime: $r(x)\ge \upkappa r(y)$, 
\item The long to short range regime: $r(y)\ge \upkappa r(x)$.
\end{enumerate}

\subsubsection{First  and second regime} In these regimes, we have $r(x)\simeq d(x,y)$ hence:
$$|R(x,y)|\le C\left[ \int_0^{r^2(x)}  \frac{ e^{-\frac{r^2(x)}{ct}}}{ \volu B(x,\sqrt{t})}\frac{dt}{t}+\int_{r^2(x)}^{+\infty}  \frac{e^{-\frac{r^2(x)}{ct}}}{r(x)\volu B(x,\sqrt{t})} \frac{dt}{\sqrt{t}}\right]\,\,.$$
Using the doubling assumption, the first integral is bounded by 
$$C \frac{r(x)^{\mu}}{\volu B(x,r(x))} \int_0^{r^2(x)} \frac{ e^{-\frac{r^2(x)}{ct}}}{t^{\frac\mu2+1}}dt\le  \frac{C}{\volu B(x,r(x))} \le  \frac{C}{\volu B(o,r(x))} \, .$$
For the second integral, we use the fact that if $\sqrt{t}\ge r(x)$ then 
$$\volu B(x,\sqrt{t}\,)\ge C \volu B(o,\sqrt{t}\,)\ge C \left(\frac{\sqrt{t}}{r(x)}\right)^\nu \volu B(o,r(x)\,)\,\,$$
and because $\nu>1$, we obtain
\begin{equation*}\begin{split}\int_{r^2(x)}^{+\infty}  \frac{e^{-\frac{r^2(x)}{ct}}}{r(x)\volu B(x,\sqrt{t})} \frac{dt}{\sqrt{t}}& \le C\frac{r^{\nu-1}(x)}{\volu B(x,r(x))}    \int_{r^2(x)}^{+\infty} \frac{ e^{-\frac{r^2(x)}{ct}}}{t^{\frac{\nu+1}{2}}}dt\\
&\le \frac{C}{\volu B(o,r(x))}.\end{split}\end{equation*}

\subsubsection{The long to short range regime} In this regime we have $r(y)\simeq d(x,y)$ hence we have
\begin{equation*}\begin{split}
|R(x,y)|\le &C \int_0^{r^2(x)}  \frac{ e^{-\frac{r^2(y)}{ct}}}{ \volu B(x,\sqrt{t})}\frac{dt}{t}\\
&+C \int_{r^2(x)}^{r^2(y)}  \frac{ e^{-\frac{r^2(y)}{ct}}}{ r(x)\volu B(x,\sqrt{t})}\frac{dt}{\sqrt{t}}+C\int_{r^2(y)}^{+\infty}  \frac{e^{-\frac{r^2(y)}{ct}}}{r(x)\volu B(x,\sqrt{t})} \frac{dt}{\sqrt{t}}\,\,,\end{split}\end{equation*}

Using the same techniques, we get 
$$\int_0^{r^2(x)}  \frac{ e^{-\frac{r^2(y)}{ct}}}{ \volu B(x,\sqrt{t})}\frac{dt}{t}\le C \frac{r(y)^{\mu}}{\volu B(x,r(y))} \int_0^{r^2(x)} \frac{ e^{-\frac{r^2(y)}{ct}}}{t^{\frac\mu2+1}}dt\le  \frac{C}{\volu B(o,r(y))}$$

and 
$$ \int_{r^2(x)}^{r^2(y)}  \frac{ e^{-\frac{r^2(y)}{ct}}}{ r(x)\volu B(x,\sqrt{t})}\frac{dt}{\sqrt{t}}\le \frac{r(y)}{r(x)} \frac{C}{\volu B(o,r(y))}\,.$$
Similarly, the reverse doubling hypothesis (RD${}_{\nu}$) and the fact that $\nu>2$ yield:
$$\int_{r^2(y)}^{+\infty}  \frac{e^{-\frac{r^2(y)}{ct}}}{r(x)\volu B(x,\sqrt{t})} \frac{dt}{\sqrt{t}}\le \frac{r(y)}{r(x)} \frac{C}{\volu B(o,r(y))}.$$

As a conclusion, we have obtained:
\begin{lem}\label{lem:riesz}
There is a positive constant $C$ such that:
\begin{itemize}
\item
if $x,y\in M$ satisfy  $d(x,y)\ge \frac{1}{\upkappa} r(x)$ and $\frac{1}{\upkappa} r(x)\le  r(y)\le \upkappa r(x)$ or  $r(x)\ge \upkappa r(y)$ then
$$
|R(x,y)|\le  \frac{C}{\volu B(o,r(x))}\,\,,
$$
\item if $x,y\in M$ satisfy $r(y)\ge \upkappa r(x)$ then
$$
|R(x,y)|\le   \frac{r(y)}{r(x)}\frac{C}{\volu B(o,r(y))}\,\,.
$$
\end{itemize}
\end{lem}

\section{boundedness of the Riesz transform}
\subsection{}In this section, the \tref{thm:general} will be proved,  hence we consider a complete Riemannian manifold $(M,g)$ satisfying the hypotheses of this theorem.

When $f\in \cC^\infty_0(M)$, we decompose $R(f)$ in three parts:
$$R(f)=R_d(f)+R_0(f)+R_1(f)\, ,$$
where the Schwartz kernels of $R_0$ and $R_1$ are locally bounded and given by the restriction of the Schwartz kernel of $R$ to the sets
$$\Omega_0:=\{(x,y)\in M\times M\,,\, d(x,y)\ge \upkappa^{-1} r(x)\,\mathrm{and}\,  \upkappa r(x)\ge r(y)\}\ ,$$
$$\Omega_1:=\{(x,y)\in M\times M\,,\,   \upkappa r(x)\le r(y)\}\ .$$
That is to say for $\alpha\in \cC^\infty_0(T^*M)$ and $f\in \cC^\infty_0(M)$:
$$\langle \alpha, R_0(f)\rangle_{L^2}=\int_{\Omega_0} \langle \alpha(x),R(x,y) \rangle_g f(y) dy dx.$$
and similarly for $R_1$. 
Note that if $\upkappa r(x)\le r(y)$ then, recalling  $\upkappa\ge 4$, we get:
$$d(x,y)\ge r(y)-r(x)\ge(\upkappa-1)r(x)\ge \upkappa^{-1} r(x)\, .$$ Hence on $\Omega_1$ we are far away from the diagonal.  The Schwartz kernel of the Riesz transform has a singularity along the diagonal of $M\times M$ and  $R_d$ is the restriction the kernel of the Riesz transform to the following neighborhood of the diagonal:
$$\cV(Diag):=\{(x,y)\in M\times M\,,\,   d(x,y)\le \upkappa^{-1} r(x)\}\ .$$
\subsection{The short to long range part}
This part is now relatively easy to handle:
\begin{prop}\label{prop:shortlong} The operator $R_0$ is bounded $L^\infty(M)\rightarrow L^\infty(T^*M)$ and $L^1\rightarrow L^1_w$: that is to say for any $f\in L^1$ and any $\lambda>0$, we have
$$\volu\{x\in M\,,\, |R_0(f)(x)|>\lambda\}\le \frac{C}{\lambda} \|f\|_{L^1}\ .$$
In particular by interpolation, $R_0\colon L^p(M)\rightarrow L^p(T^*M)$ is bounded for any $p\in (1,+\infty)$.
\end{prop}
\proof
As a matter of fact, our previous analysis (\lref{lem:riesz}) implies that if $f\in \cC^\infty_0(M)$ then
$$ |R_0(f)(x)|\le \frac{C}{\volu B(o,r(x))} \int_{B(o,\upkappa r(x))}| f(y)|dy\ ,$$
hence the boundedness  $L^\infty(M)\rightarrow L^\infty(T^*M)$ is a direct consequence of the doubling property.
Moreover this also implies that 
$$\{x\in M\,,\, |R_0(f)(x)|>\lambda \}\subset B(o,\rho)$$ where
$\rho$ satisfies$$\volu B(o,\rho)< \frac{C}{\lambda} \|f\|_{L^1}\, $$
and the boundedness $L^1\rightarrow L^1_w$ follows immediately.
\endproof
\subsection{The diagonal part} In this part, we are going to use an idea from \cite[section 4]{CD} and a result from \cite[section 4]{ACDH} in order to prove the following:
\begin{prop}\label{AC} The operator $R_d$ is bounded on $L^p$ for every $p\in (2,+\infty)$.
\end{prop}
\noindent\textit{Proof of \pref{AC}.}
We build a cover on $M$  by remote balls. By induction on $N\in \N$:
\begin{itemize}
\item $B_{0,1}=B(o,1)$.
\item We cover $B(o,2^N)\setminus \cup_{i<N, j} B_{i,j}$ by a collection of  balls $B_{N,1},\dots ,B_{N,k_N}$ of radius $2^{N-10}$ that are centered on  $B(o,2^N)\setminus B(o,2^{N-1})$ and such that the balls $\frac12 B_{N,1},\dots ,\frac12 B_{N,k_N}$ are disjoint and included in $B(o,2^N)\setminus \cup_{i<N, j} B_{i,j}$. 
\end{itemize} 
At each stage $N$, the number of balls is bounded independently of $N$: $$k_N\le m(\upvartheta).$$
We obtain in this way a subset  $A\subset \N^2$ and a cover $$M=\cup_{\alpha\in A} B_\alpha$$  by balls 
$B_\alpha=B(x_\alpha,r_\alpha)$. Note that we have by construction:
$2^{-10} r(x_\alpha)\le r_\alpha\le 2^{-9} r(x_\alpha).$ Moreover this cover has a finite multiplicity: there is a constant $p$ such that for any $x\in M$: 
$$\card\{\alpha\in A, x\in B_\alpha\}\le p.$$

Let $\chi_\alpha$ be a partition of unity subordinate to this covering; we can assume $|d\chi_\alpha|\le C/r_\alpha$.
If $\upkappa$ is chosen large enough $(\upkappa\ge 2^{10})$, then we have
$$|R_d(f)(x)|\le \sum_\alpha |\un_{4B_\alpha}(x)\, R(\chi_\alpha f)(x)|.$$
Let $R_\alpha=\un_{4B_\alpha}R\chi_\alpha$, we decompose 
$$R_\alpha=R_{\alpha,0}+R_{\alpha,1}\ ,$$
where 
$$R_{\alpha,0}(f)(x)=\un_{4B_\alpha}(x)\int_0^{r_\alpha^2} \nabla_xe^{-\tau \Delta}(\chi_\alpha f)(x)\frac{d\tau}{\sqrt{\pi \tau}}$$ and
$$R_{\alpha,1}(f)(x)=\un_{4B_\alpha}(x)\int_{r_\alpha^2}^\infty \nabla_xe^{-\tau \Delta}(\chi_\alpha f)(x)\frac{d\tau}{\sqrt{\pi \tau}}\, .$$ 
The covering $M=\cup_{\alpha\in A} B_\alpha$ has finite multiplicity and  $(M,g)$ is doubling, hence we 
only need to prove that there is a uniform constant $C$ such that for all $\alpha$
$$\|R_{\alpha,0}\|_{L^p\to L^p}\le C\,\,,\,\, \|R_{\alpha,1}\|_{L^p\to L^p}\le C.$$

\begin{lem}\label{R0}
There is a constant $C$ independent of $\alpha$ such that:$$\|R_{\alpha,0}\|_{L^p\to L^p}\le C.$$
\end{lem}
\proof We will use the arguments of \cite[subsection 3.2]{ACDH} and \cite[section 4]{ACDH} together with 
 the following estimates on the gradient of the heat kernel: for all $x,y\in M$ and all  $t\in (0,Ar^2(x))$
\begin{equation}\label{gradientheat}|\nabla_xh_t(x,y)|\le \frac{C}{\sqrt{t}\, \volu B(x,\sqrt{t})} e^{-\frac{d^2(x,y)}{ct}}.\end{equation}

We will apply  \cite[Theorem 2.4]{ACDH}. The setting is the following:

\begin{itemize}
\item $(M,g)$ is a complete Riemannian manifold,
\item $T\colon L^2(M)\rightarrow L^2(M)$ is a bounded sublinear operator, 
\item $\{A_r\}_{r>0}$ is a family of bounded operator on $L^2$:
$$\sup_{r>0} \|A_r\|_{L^2\to L^2}<\infty.$$
\item $U\subset \Omega\subset M$ are two open subsets such that $\Omega$ satisfies the relative doubling condition:
there is a constant $\bar\upvartheta $ such that for all balls $B\subset M$:
$$\volu (2B\cap \Omega)\le \bar\upvartheta  \volu(B\cap \Omega).$$
\item $S\colon L^p(U)\rightarrow L^p(\Omega)$ is a bounded operator for all $p>2$.
\end{itemize}

The assumption are 
\begin{enumerate}[i)]
\item For all $p>2$, the  sublinear operator $\cM^\#$ defined by 
$$\cM^\#(f)(x)=\sup_{B, B\cap \Omega\ni x}\frac{1}{\volu (\Omega\cap B)}\int_{B\cap \Omega} |T(\Id-A_{r(B)})(f)|^2$$
is bounded on $L^p$.
\item For all $f\in L^p(U)$ and all balls $B\subset M$ and $x,y\in \Omega\cap B$
$$|TA_{r(B)}(f)|^2(y)\le C \cM_\Omega\left(|T(f)|^2+|S(f)|^2\right)(x).$$
\end{enumerate}
 Where we have noted $\cM_\Omega$ the maximal operator relative to $\Omega$:
 $$\cM_\Omega(f)(x)=\sup_{B, B\cap\Omega\ni x}\frac{1}{\volu (\Omega\cap B)}\int_{B\cap \Omega} |f|\, .$$
 
 The conclusion is  that the operator $$T\colon L^p(U)\rightarrow L^p(\Omega)$$ is  bounded and there is an upper bound on operator norm of $T\colon L^p\to L^p$ that depends only on the constants involved in the setting and the hypothesis.
 
The following result can be deduced from \cite[Proposition 2.4]{DR} or \cite[p 1159]{CD}: 
 
 \begin{lem}\label{firstlemm}If $(M,g)$ is doubling then
 $$\int_{M\setminus B(x,r)} e^{-\frac{d^2(x,y)}{c t}} |f(y)|dy\le C \volu B(x,\sqrt{t})\, e^{-\frac{r^2}{2ct}}\, \cM(f)(x)\, ,$$

 where $$\cM(f)(x)=\sup_{B, B\ni x}\frac{1}{\volu ( B)}\int_{B} |f|$$ is the maximal operator associated to $(M,g)$.
  \end{lem}
 The next lemma is folklore and a proof can be found in \cite[Proposition 4.14]{Cheediff} or \cite[lemma 2.3.1]{KZ}:

  \begin{lem}\label{secondlemm}I If $B(x,r)\subset M$ is a remote ball and if $M$ satisfies the (QD) condition, then for all $f\in \cC^1(B(x,r)):$
  $$\left| f(x)-f_{B(x,r)}\right|\le C \, r\, \sup_{0<s\le r} \frac{1}{\volu ( B(x,s))}\int_{B(x,s)} |df|. $$
  \end{lem}
 
  Using the lemma \ref{firstlemm}, the  gradient estimate (\ref{gradientheat}) and the argumentation of \cite[subsection 3.2]{ACDH} we easily get:
  $$\cM^\#(f)(x)\le C\sqrt{ \cM_\Omega( |f|^2)(x)\,\,}.$$
  
Moreover if $f\in \cC^1_0(M)$, if $B$ is a remote ball of radius $r$ and if $x,y\in B$ then the lemma \ref{secondlemm}, yields
\begin{equation}\label{Maxine}
\left|\nabla e^{-r^2\Delta} f\right |(y)\le C\sup_{0<s\le r} \frac{1}{\volu ( B(x,s))}\int_{B(x,s)} |df|+\frac{C}{r}\cM(|f|)(x).\end{equation}
  
  We will use \cite[Theorem 2.4]{ACDH} with $U=B_\alpha$, $\Omega=4B_\alpha$, $A_r=e^{-r^2\Delta}$ and 
  $$T(f)(x)=\left|\int_0^{r_\alpha^2} \left(de^{-\tau \Delta}f\right)(x) \frac{d\tau}{\sqrt{\pi\tau\,}}\right| .$$
  Let $v\in L^p(U)$, if we apply the last inequality (\ref{Maxine}) to $$f=\cS(v)=\int_0^{r_\alpha^2} e^{-\tau \Delta}v \frac{d\tau}{\sqrt{\pi\tau\,}}\,\,,$$
  then we get that for all  balls  $B$ and all $x,y\in B\cap\Omega$:
  $$|TA_{r(B)}v|(y)\le C\sup_{0<s\le r_\alpha} \frac{1}{\volu ( B(x,s))}\int_{B(x,s)} |T(v)|+\frac{C}{r_\alpha}\cM(\cS(v))(x).$$
  The estimate of the gradient of the heat kernel implies that 
  $$\int_{B(x,s)} |\cS(v)|+\int_{B(x,s)} |T(v)|\le\int_{B(x,s)\cap \Omega} |T(v)|+C\frac{ \volu (B(x,s)\setminus \Omega)}{\volu U}\int_U |v|\,\,.$$
  Hence we get 
  $$|TA_{r(B)}v|(y)\le C \cM_\Omega(T(v))(x)+\frac{C}{r_\alpha} \cM_{\Omega}(\cS(v))(x)+c\cM_{\Omega}(v)(x)\, .$$
  The fact that 
  $$\|\cS\|_{L^p\to L^p}\le C\, r_\alpha$$ yields the uniform estimate
  $$\|R_{\alpha,0}\|_{L^p\to L^p}\le\|T\|_{L^p(U)\to L^p(\Omega)}\le C.$$
  \endproof
  
  \begin{lem}\label{R1}There is a constant $C$ independent of $\alpha$ such that:$$\|R_{\alpha,1}\|_{L^p\to L^p}\le C.$$
  \end{lem}
  
\proof We use the fact that for all $y\in B_\alpha$, $x\in 4B_\alpha$  and $t\ge r_\alpha^2$ then
$$|\nabla_xh_t(x,y)|\le \frac{C}{r_\alpha\, \volu B(x,\sqrt{t})} e^{-\frac{d^2(x,y)}{ct}}.$$
So that the Schwartz kernel of $R_{\alpha,1}$ is bounded by
$$C \un_{4B_\alpha}(x)\un_{B_\alpha}(y)\int_{r_\alpha}^\infty\frac{e^{-\frac{d^2(x,y)}{ct}}}{r_\alpha
 \volu B(x,\sqrt{t})} e^{-\frac{d^2(x,y)}{ct}} \frac{dt}{\sqrt{t}}\, .$$
 And using the slow variation of the volume of balls, we get that for all $x\in 4B_\alpha$ and all   $t\ge r^2_\alpha$:  $\volu B(x,\sqrt{t})\simeq \volu B(o,\sqrt{t})$. And with the reverse doubling condition (RD${}_{\nu}$), we obtain:
 $$|T_{\alpha,1}(f)(x)|\le \frac{C\un_{4B_\alpha}(x)}{\volu B_\alpha}\int_{B_\alpha} |f|(y)dy.$$
 Hence $T_{\alpha,1}$ is bounded on $L^1$ and on $L^\infty$ with an operator norm bounded independently of $\alpha$.  \endproof
 
\noindent The \lref{R0} and  \lref{R1} imply the \pref{AC}\hfill$\square$
\subsection{The long to short range part} This is the most significant part ; according to our previous analysis,
the $L^p$ boundedness of the Riesz transform is equivalent to the $L^p$ boundedness of $R_1$. And with the \lref{lem:riesz}, it is sufficient  to find conditions under which the operator
$$T(f)(x):=\frac{1}{r(x)}\int_{M\setminus B(o,\upkappa r(x))} \frac{r(y)}{\volu B(o,r(y))} |f(y)| dy$$ is bounded $L^p\rightarrow L^p_w$.

 When $f\in L^p(M)$, we have:
$$|T(f)(x)|\le M_p(x)\|f\|_{L^p}$$
where 
$$M_p(x)=\frac{1}{r(x)}\left(\int_{M\setminus B(o,\upkappa r(x))} \left(\frac{r(y)}{\volu B(o,r(y))} \right)^{\frac{p}{p-1}}dy\right)^{1-\frac1p}.$$
If we introduce the Riemann-Stieljes measure associated to the non-decreasing function $V(r)=\volu B(o,r)$, we get (by integrating by parts)
\begin{equation*}\begin{split}
\int_{M\setminus B(o, R)} \left(\frac{r(y)}{\volu B(o,r(y))} \right)^{\frac{p}{p-1}}dy&=\int_R^\infty \left(\frac{r}{\volu B(o,r)} \right)^{\frac{p}{p-1}}dV(r)\\
&=(p-1)\frac{R^{\frac{p}{p-1}}}{V(R)^{\frac{1}{p-1}}}+p\int_R^\infty \frac{r^{\frac{1}{p-1}}}{V(r)^{\frac{1}{p-1}}}dr.
\end{split}\end{equation*}
provided \begin{itemize}
\item $\displaystyle \lim_{R\to \infty}\frac{R^{\frac{p}{p-1}}}{V(R)^{\frac{1}{p-1}}}=0$ and
\item $\displaystyle \int_1^\infty\left( \frac{r}{V(r)}\right)^{\frac{1}{p-1}}dr<\infty$.
\end{itemize}
The second condition implies the first one and 
in our setting, the second condition is equivalent to the $p$-hyperbolicity of the manifold $(M,g)$ (\cite{HK}).  Recall that  if $\cO\subset M$ then its $p$-capacity is defined by: 
$$\capa_p \cO:=\inf\left\{\int_M |d\varphi|^p\vol_g\,,\, \varphi\in \cC^\infty_0(M) \,\, \mathrm{and}\,\, \varphi\ge 1\,\, \mathrm{on}\,\, \cO\right\}\,\,,$$ and that  $(M,g)$ is said to be $p$-hyperbolic if the $p$-capacity of bounded open subsets is positive. 
Using the argument of the proof of the \pref{prop:shortlong}, we obtain as before that $$T\colon L^p\rightarrow L^p_w$$ is bounded provided that for some constant $C$ independent of $R$, we have:
\begin{equation}\label{condip}
 \int_R^\infty\left( \frac{r}{V(r)}\right)^{\frac{1}{p-1}}dr\le C \left(\frac{r^p}{V(r)}\right)^{\frac{1}{p-1}}.
\end{equation}
 
Using the reverse doubling condition (RD${}_{\nu}$), we get that this condition is satisfied when $p<\nu$; hence we have proven the \tref{thm:general}.

According to \cite{CSCH}, in our setting the $p$-capacity of an anchored ball can be estimated:
$$\capa_p (B(o,R))\le C \left(\int_R^\infty \frac{r^{\frac{1}{p-1}}}{V(r)^{\frac{1}{p-1}}}dr\right)^{1-p}\,\,.$$
Hence, when the $p$-capacity of anchored balls satisfy the uniform estimate:
$$\capa_p(B(o,R))\ge c \frac{\volu B(o,R)}{R^p}$$ then 
the condition (\ref{condip}) is satisfied. This argumentation provides a direct proof of  the \cref{coro:general}.

\section{Passing the volume growth exponent}
We now are proving  the \tref{thm:impro}.
\subsection{Reverse H\"older inequality}

According to \cite{GS0}, when $(M,g)$ satisfies the conditions (D), (QD) and (RCA) then $(M,g)$ satisfies the  scale invariant $L^1$ Poincar\'e inequalities: for any ball $B\subset M$ and any $f\in \cC^\infty(2B)$:
$$\|f-f_B\|_{L^1(B)} \le C r(B)\, \|df\|_{L^1(2B)}.$$

These Poincar\'e inequalities and the doubling condition implies (\cite{AC}) that for all $q\in (1,2]$, the reverse Riesz transform is bounded in $L^q$:
there is a constant $C$ such that for any $f\in \cC^\infty_0(M)$:
$$\|\sqrt{\Delta} f \|_{L^q}\le C \|df\|_{L^q}.$$
Hence when $p>2$, the Riesz transform is bounded on $L^p$ 
as soon as the Hodge projector $\Pi=d\Delta^{-1} d^*\colon L^2(T^*M)\rightarrow  L^2(T^*M)$
 has a bounded extension to $L^p$ (cf. \cite[lemma 0.1]{AC}).

Now the proof of the implication $1)\Rightarrow 2)$ of the \cite[theorem 2.1]{AC} (see also \cite{shen}) shows that if for some $\tilde p> p$,  we have a $L^{\tilde p}$-reverse H\"older inequality for the gradient of harmonic functions, then the 
Hodge projector has a bounded extension on $L^p$.
\begin{defi}
 A complete Riemannian manifold $(M,g)$ is said to satisfy the $L^p$-reverse H\"older inequality if for 
some constants $C>0$, $\bar\alpha>\underline{\alpha}>1$ and for any ball $B\subset M$ and any harmonic function $h$ defined on $\bar\alpha B$, one has the reverse H\"older inequality:
$$\left(\frac{1}{\volu B}\int_B |dh|^p\right)^{\frac{1}{p}} \le
 C \left(\frac{1}{\volu(\underline{\alpha}B)}\int_{\underline{\alpha}B} |dh|^2\right)^{\frac{1}{2}} .$$
\end{defi}

In our case, the quadratic decay  of the negative part of the Ricci curvature and the Cheng and Yau's estimate on the gradient of harmonic function (\cite{ChengYau}) yield a $L^\infty$-reverse H\"older inequality for remote balls. The following lemma shows that in our setting, we will get  the $L^{ p}$-reverse H\"older inequality provided it holds for anchored balls.
\begin{lem}\label{lem:anchoresu}
Let  $(M,g)$ be a complete Riemannian manifold that satisfies the doubling condition. 
The $L^{ p}$-reverse H\"older inequality holds provided it holds for remote and anchored balls.
\end{lem}
 \proof Assume that the $L^{ p}$-reverse H\"older inequality holds for remote and anchored balls with parameters $\bar\alpha>\underline{\alpha}>1$.

Let $B(x,r)$ be a ball that is not anchored nor remote, i.e. $x\not= o$ and $r\ge r(x)/2$. Let $\uplambda\ge 1$ be a real parameter.
\begin{enumerate}[i)]
 \item Assume that $r\ge \uplambda r(x)$ and let $B'=B(o,(1+\uplambda^{-1}) r)$
we get 
$B(x,r)\subset B'$ and
$\alpha B'\subset \beta B(x,r)$ provided $\beta=(1+\uplambda^{-1}) \alpha+ \uplambda^{-1}$.
Define now 
  $\bar\beta= (1+\uplambda^{-1}) \bar \alpha+ \uplambda^{-1}$ and $\underline{\beta}=(1+\uplambda^{-1}) \underline{\alpha}+ \uplambda^{-1}$.
The six balls $$B(x,r)\,\,,\,\, \underline{\beta}B(x,r)\,\,,\,\, \bar\beta B(x,r)\,\,,\,\,
 B'\,\,,\,\,\underline{\alpha}B'\,\,,\,\,\bar \alpha B'$$ have a comparable volume. The inclusions
$B(x,r)\subset B'$, $\bar\alpha B'\subset \bar\beta B(x,r)$ and $\underline{\alpha} B'\subset \underline{\beta} B(x,r)$ together with $L^{ p}$-reverse H\"older inequality for the ball $B'$ imply 
that if $h$ is a harmonic function defined on $\bar\beta B(x,r)$, then 
$$\left(\frac{1}{\volu B(x,r)}\int_{B(x,r)} |dh|^p\right)^{\frac{1}{p}} \le
 C \left(\frac{1}{\volu(\underline{\beta}B(x,r))}\int_{\underline{\beta}B(x,r))} |dh|^2\right)^{\frac{1}{2}} .$$
\item Assume now that $\frac{r(x)}{2}\le r\le  \uplambda r(x)$: 
let $h$ be a harmonic function defined on $B(x,4r)$. Then we have the inclusion $B(o,(4- \uplambda^{-1})r)\subset B(x,4r)$.
We consider a minimal covering of  $B(o,(4- \uplambda^{-1})r)\setminus B(o,4\delta r)$ by balls of radius $\delta r$:
$$B(o,(4- \uplambda^{-1})r)\setminus B(o,4\delta r)=\cup_{i\in I} B_i\, .$$
All the balls $B_i$ are remote and for some constant $N$ depending only on $\delta$ and on the doubling constant $\upvartheta$ 
$$\card I\le N.$$
Moreover all the balls $B_i$, $B(o,4\delta r)$ have a comparable volume.
We  choose $\delta$ so that 
$$ 8\delta<1\,\,\mathrm{and}\,\,(1+\bar\alpha) \delta<1  \, .$$

We introduce the collection of balls: $\cB_*=\left\{B(o,4\delta r)\right\}\cup \left\{B_i,\ i\in I\right\}$ and let 
$$\cB=\{B\in \cB_* , B\cap B(x,r)\not=\emptyset\}\,\,.$$
If $B\in \cB$ then $\underline{\alpha}B\subset B(x,2r)$
and also $\bar\alpha B\subset B(x,2r)$, so that: 
\begin{equation*}\begin{split}\int_{B(x,r)}|dh|^p &\le \sum_{B\in \cB} \int_{B}|dh|^p\\
&\le C \sum_{B\in \cB}  (\volu B)^{1-\frac 2p} \left(\int_{\underline{\alpha} B}|dh|^2\right)^{\frac2p}\\
&\le C\, (\card  I+1)(\volu B(x,r))^{1-\frac 2p} \left(\int_{2 B(x,r)}|dh|^2\right)^{\frac2p}\,\,.
\end{split}\end{equation*}
\noindent Hence the result.
\end{enumerate}

 \endproof
\subsection{H\"older Elliptic estimates and the  Green kernel }

Recall that we say that $(M^n,g)$ satisfies the scale invariant H\"older Elliptic estimates (EH${}_\alpha$) if there is a constant $C$ such that for any ball $B\subset M$ and any harmonic function $h$ defined on $3B$,  then we have for all $x,y\in B$:
$$|h(x)-h(y)|\le C \,\left(\frac{d(x,y)}{r(B)}\right)^\alpha\, \sup_{z\in 2B} |h(z)|\,\,.$$

\begin{rem} The argument  given in the proof of \lref{lem:anchoresu}) shows that with the (RCA) and the doubling condition, the  scale invariant H\"older Elliptic estimates (EH${}_\alpha$) holds for all balls provided its holds for all remote and anchored balls. In our setting, the quadratic decay  of the negative part of the Ricci curvature implies a Lipschitz  estimates for harmonic function on remote balls;  hence the scale invariant $\alpha$-H\"older Elliptic estimates holds for any ball if and only if it holds for anchored balls.\end{rem}

We assume now that $(M,g)$ is a complete Riemannian manifold that satisfies the conditions (QD), (VC) and (RCA) and the scale invariant H\"older Elliptic estimates (EH${}_\alpha$).
For $R\ge 1$ and $\upkappa> 2$, we consider the anchored balls 
$$B^*=B(o,R/\upkappa)\subset \subset B^\#=B(o, R). $$

We will note $h^\#(t,x,y)$ the heat kernel on the ball $B^\#$ for the Dirichlet boundary condition and $G^\#(x,y)$ the associated Green kernel:
$$G^\#(x,y)=\int_0^{+\infty} h^\#_t(x,y)dt\,\,. $$
If $f\in L^2(B^\#)$, then 
$$u(x)=\int_{B^\#}G^\#(x,y)f(y)dy$$ is the solution of the equation:
$$\Delta u=f\hspace{0,2cm}:\hspace{0,2cm} u=0\hspace{0,2cm} \mathrm{on}\hspace{0,2cm} \partial B^\#\,.$$

\noindent Because our manifold satisfies the Faber-Krahn inequality, we get 
$$\lambda^D_1(B^\#)\ge \frac{c}{R^2}\,.$$
Recall that for all $t>0, x,y\in B^\#$:
$$h^\#_t(x,y)\le h_t(x,y),$$
and for all $t>s>0, x,y\in B^\#$:
$$h^\#_t(x,y)\le \sqrt{h^\#_t(x,x)\,h^\#_t(y,y)\,}\le e^{-\lambda_1^D(B^\#)(t-s)} \sqrt{h^\#_s(x,x)\,h^\#_s(y,y)\,}\,.$$

\noindent Hence we get the estimate:
$t\in (0,R^2), x,y\in B^\#$:
$$h^\#_t(x,y)\le \frac{C}{\volu B(x,\sqrt{t})} e^{-\frac{d^2(x,y)}{ct}}\, .$$
 And for $t\ge R^2, x,y\in B^\#$:
$$h^\#_t(x,y)\le \frac{C}{\volu B^\#} e^{-c\frac{t}{R^2}}\,.$$

\noindent Hence if $x\in B^*, y\in B^\#\setminus \frac12B^\#$ we get $d(x,y)\simeq R$ and 
$$G^\#(x,y)\le \int_0^{R^2} \frac{C}{\volu B(x,\sqrt{t})} e^{-\frac{R^2}{ct}}dt+\frac{CR^2}{\volu B^\#}\,\,.$$
But  using the doubling condition we obtain: 
$$\int_0^{R^2} \frac{C}{\volu B(x,\sqrt{t})} e^{-\frac{R^2}{ct}}dt\le \int_0^{R^2} \frac{C}{\volu B(x,R)}\left(\frac{R}{\sqrt{t}}\right)^\mu e^{-\frac{R^2}{ct}}dt \le  \frac{CR^2}{\volu B^\#}\,.$$
Finally, we obtain that if $x\in B^*$ and $y\in B^\#\setminus \frac12 B^\#$ then 
$$G^\#(x,y)\le C\frac{R^2}{\volu B^\#}\,\,.$$

Using the  $\alpha$-H\"older regularity estimate, we get the following estimation on the gradient of the Green kernel:

\begin{lem}\label{lem:impr}  There is a   constant $C>0$ such that if $x\in \frac12 B^*$ and $y\in B^\#\setminus \frac12 B^\#$ then
$$|\nabla_xG^\#(x,y)|\le C \left(\frac{R}{r(x)+1}\right)^{1-\alpha}\frac{R}{\volu B^\#}.$$
\end{lem}
\proof Indeed we apply the above $\alpha$-H\"older regularity estimates for the harmonic function $h(z)=G^\#(z,y)-G^\#(x,y)$. We have seen that 
$$\sup_{z\in B^*} |h(z)|\le C\frac{R^2}{\volu B^\#}\,.$$
\noindent Hence if $x\in \frac12 B^*$ and $z\in B(x,\frac12 r(x))$ we obtain
$$ |h(z)|= |h(z)-h(x)|\le C \left(\frac{r(x)}{R}\right)^\alpha \frac{R^2}{\volu B^\#}.$$
And with the Cheng and Yau's gradient estimate \cite{ChengYau}, we get:
$$|\nabla h|(x)=|\nabla_x G^\#(x,y)|\le   \frac{C}{r(x)}  \left(\frac{r(x)}{R}\right)^\alpha \frac{R^2}{\volu B^\#}.$$
When $r(x)\le 1/2$, the same idea leads to the estimates:
$$ \forall z\in B(o,1),  |h(z)|= |h(z)-h(o)|\le C \left(\frac{1}{R}\right)^\alpha \frac{R^2}{\volu B^\#}.$$
The Cheng and Yau gradient estimate implies that 
$$x\in B(o,1/2)\colon |\nabla_x G^\#(x,y)|\le  C  \left(\frac{1}{R}\right)^\alpha \frac{R^2}{\volu B^\#}.$$
\endproof

\subsection{From  $\alpha$-H\"older elliptic regularity to reverse H\"older inequality}
In order to prove the \tref{thm:impro}, we only need to prove that 
if $(M,g)$ satisfies the hypotheses of this theorem and if $p$ is such that $(1-\alpha)p<\nu$ then  anchored balls satisfy the $L^{ p}$-reverse H\"older inequality.

With J.Cheeger and T.Colding \cite{ChCo}, we can build a smooth function $\chi$ such that
 \begin{itemize}
\item $\chi=0$ on $M\setminus B(o,3R/4)$
\item $\chi=1$ on $ B(o,R/2)$ 
\item $R|d\chi|+R^2|\Delta \chi|\le C$.
\end{itemize}

Indeed, according to \cite[Theorem 6.33]{ChCo}, if $B$ is a remote ball then there is a smooth function $\chi_B$ with compact support in $B$ such that $\chi_B=1$ on $\frac12 B$ and such that
$$ r(B)|d\chi_B|+r^2(B)|\Delta \chi_B|\le C(n,\kappa).$$
Using the doubling hypothesis, we can cover $\partial B(o,R/2)$ by at most $N$ balls of radius $R/8$ and centered in $\partial B(o,R/2)$:
$$\partial B(o,R/2)\subset \cup_i B_i\, .$$
Introduce $\varphi=\sum_i \chi_{2B_i}$ on $\partial B(o,R/2)$, we have $\varphi\ge 1$ and $\varphi$ has compact support in $B(o,\frac34R)$.
Moreover there is a constant (independent of $R$) such that:
$R|d\varphi|+R^2|\Delta \varphi|\le C$.
We let $\chi=u(\varphi)$ where $u\colon [0,\infty[\rightarrow [0,1]$ is a smooth function such that $u=1$ on $[1,\infty)$ and $u=0$ on $[0,1/2]$.

We will use that annuli satisfy a scale invariant Poincar\'e inequality (according to \cite{GS0}):
 there are constants $\upkappa>1$ and $C>0$ such that for any $R>1$  then  
$$\forall f\in \cC^\infty(A_R^*) \colon \|f-f_{A_R}\|_{L^2(A_R)}\le C R  \|df\|_{L^2(A_R^*)}\, ,$$
where $$A_R:=B(o, R)\setminus B(o,R/2) \ \mathrm{and}\  
A_R^*:=B(o,\upkappa R)\setminus B(o,R/(2\upkappa)).$$
Let $B^\#=B(o, R)$ and let $G^\#(x,y)$ be the Green kernel of the Dirichlet Laplacian on $B^\#$.

Let $h$ be a harmonic function on $B(o,\upkappa R)$.
 Using the Green kernel, we have for 
$x\in B(o,R/(2\upkappa) )$  and any constant $c$:
$$h(x)-c=\int_{A_R} G^\#(x,y)\Delta(\chi (h-c))(y)dy\, .$$
Using the fact that $h$ is harmonic, we know that 
$$ \Delta(\chi (h-c))=(\Delta\chi)(h-c)-2\langle d\chi,dh\rangle_g\,\,.$$
With Cauchy-Schwarz inequality, we obtain
\begin{equation*}\begin{split}
  |dh|^2(x) \le &\,C\int_{A_R} |\nabla_xG^\#(x,y)|^2 dy\\
&\times 
 \int_{A_R} \left[R^{-4}|h(y)-c|^2+R^{-2} |dh|(y)^2\right] dy \,\,.
\end{split}\end{equation*}
Choosing $$c=\frac{1}{\volu(A_R) }\int_{A_R}h(y)dy=h_{A_R}$$ and
using the Poincar\'e inequality for the annulus, we get 
$$ \int_{A_R}|h(y)-c|^2dy \le  C R^{2}\int_{A_R^*}|dh|(y)^2dy\,\,.$$
Eventually, we obtain 
$$|dh|^2(x) \le C R^{-2} \left(\int_{A_R} |\nabla_xG^\#(x,y)|^2 dy \right) \|dh\|_{L^2(B(o,\upkappa R))}^2\,\,.$$
Recall the estimates of the gradient of the Green kernel: for $x\in B(o,R/(2\upkappa))$ and $y\in B(o,R)\setminus B(o,R/2)$, we have
$$|\nabla_xG^\#(x,y)| \le C \left(\frac{R}{r(x)+1}\right)^{1-\alpha}\,\frac{R}{\volu B(o,R)}\,\,.$$
So that we get that 
$$ \left(\frac{1}{\volu B(o,R/(2\upkappa))}\int_{B(o,R/(2\upkappa))} |dh|^p\right)^{\frac{1}{p}} \le
 \left(\frac{ C(R)^2}{\volu(B(o,\upkappa R))}\int_{B(o,\upkappa R))} |dh|^2\right)^{\frac{1}{2}}, $$
where 
$$C^p(R)=\frac{1}{\volu B(o,R/(2\upkappa))}\int_{B(o,R/(2\upkappa))} \left(\frac{R}{r(x)+1}\right)^{p(1-\alpha)}\vol(x).$$

Now it is easy to show that the reverse doubling assumption (RD${}_{\nu}$) yields an uniform bound on $C(R)$ as soon as 
$$p(1-\alpha)<\nu.$$
\subsection{On the H\"older elliptic regularity estimates}
In our setting, the scale invariant $\alpha$-H\"older Elliptic regularity estimates are equivalent to a "quasi"-monotonicity result for the $L^2$
norm of the gradient of harmonic function.
\begin{prop}\label{prop:freq}
Assume that $(M,g)$ is a complete Riemannian manifold satisfying the hypothesis (QD), (VC) and (RCA). Let $\alpha\in (0,1]$.  Then $(M,g)$ satisfies the $\alpha$-H\"older Elliptic regularity estimates if and only if there are some constants $\upkappa>1$, $C>0$ 
 such that for any harmonic function $h$ on $B(o,\upkappa R)$ and any $1\le r\le R $
$$\frac{r^{2-2\alpha}}{\volu B(o,r) }\int_{B(o,r)} |dh|^2\le C \frac{R^{2-2\alpha}}{\volu B(o,R) }\int_{B(o,R)} |dh|^2.$$
\end{prop}
\proof If $f$ is a continuous function on a subset $\cO\subset M$, we define its oscillation by:
$$\Osc_{\cO} f=\sup_{x\in \cO} f(x)-\inf_{y\in\cO} f(y).$$
Let $R\ge 1$ and $h$ be a harmonic function defined on $B(o,2\upkappa R)$.
We let again $A_R:=B(o, R)\setminus B(o,R/2)$ and 
$A_R^*:=B(o,\upkappa R)\setminus B\left(o,R/(2\upkappa)\right)$.
We have the Poincar\'e inequality
$$\|h-h_{A_R}\|^2_{L^2(A_R)}\le C R^2 \|dh\|^2_{L^2(A_R^*)}
\, .$$
But according to \cite[lemma 6.3]{GS0}, we have  a Harnack inequality on annuli, so that 
\begin{equation}
 \label{ine1}(\Osc_{A_R} h)^2 \le \frac{C}{\volu A_R}\|h-h_{A_R}\|^2_{L^2(A_R)}\le C\frac{R^2}{\volu B(o,\upkappa R)} \|dh\|^2_{L^2(B(o,\upkappa R)}.
 \end{equation}

Using the function  $\chi$  defined by 
$$\chi(x)=\begin{cases}
1&\,\,\mathrm{on}\,\, B(o,R/2)\\
2-\frac{2r(x)}{R}  &\,\,\mathrm{on}\,\, B(o,R)\setminus B(o,R/2)\\
0&\,\,\mathrm{on}\,\, M\setminus B(o, R)\,,
\end{cases}$$
 we get 
\begin{equation*}\begin{split}\int_{B(o,R/2)} |dh|^2&\le \int_{B(o, R)} |d(\chi (h-h_{A_R}))|^2\\
&=\int_{A_R} (h-h_{A_R})^2\,|d\chi|^2 =\frac{4}{R^2}\int_{A_R} (h-h_{A_R})^2.
\end{split}\end{equation*}
In particular we get 
\begin{equation}
\label{ine2}\frac{1}{\volu B(o,R/2)} \int_{B(o,R/2)} |dh|^2\le \frac{C}{R^2} (\Osc_{A_R} h)^2\,.\end{equation}

But in our case, the $\alpha$-H\"older Elliptic regularity estimates is equivalent to
 a monotonicity inequality for $\rho\mapsto \rho^{-\alpha} \Osc_{A_\rho} h$. The result is then a consequence of the inequalities (\ref{ine1}) and (\ref{ine2}).
\endproof
\section{Non boundedness of the Riesz transform}In this section, we give two criterions for the $L^p$ unboundedness of the Riesz transform. 

\subsection{Parabolicity and the Riesz transform} Our argument is a slight improvement of some earlier results proved in \cite{CCH}. The starting point is to understand the $L^p$ closure
of the space of differential of smooth function with compact support; the following lemma is a $L^p$ adaptation of an idea from \cite{carronjga}.
\begin{lem}\label{lem:primi}
 Let $(M^n,g)$ be a complete Riemannian manifold and let $f\in W^{1,p}_{loc}$ such that $df\in L^p$.
If $M(r)=\int_{B(o,r)} |f|^p$ satisfies 
$$\int_1^\infty  \left(\frac{r}{M(r)}\right)^{\frac{1}{p-1}}dr=\infty\,,$$
then there is a sequence $(\chi_\ell)$ of smooth function with compact support such that 
$$\lim_{\ell\to\infty}\|df-d(\chi_\ell f)\|_{L^p}=0\,$$
in particular:
$$df\in \overline{d\cC^\infty_0(M)}^{L^p}.$$
\end{lem}
\proof
Let $r<R$, we define a function $\chi_{r,R}$ by
letting $\chi_{r,R}=1$ on $B(o,r)$, $\chi_{r,R}=0$ outside $B(o,R)$ and for $x\in B(o,R)\setminus B(o,r)$: 
$$\chi_{r,R}(x)=\xi_{r,R}(r(x))=\varepsilon(r,R) \int_{r(x)}^R \left(\frac{s}{M(s)}\right)^{\frac{1}{p-1}}ds\,\,, $$
where $$\varepsilon(r,R) =\left(\int_{r}^R \left(\frac{s}{M(s)}\right)^{\frac{1}{p-1}}ds\right)^{-1}.$$
Let $f$ be a function satisfying the hypotheses of the lemma, then we get:
$$\|df-d(\chi_{r,R} f)\|^p_{L^p}\le C\left( \int_{M\setminus B(o,r)} |df|^p+
\int_{B(o,R)\setminus B(o,r)} |f|^p|d \chi_{r,R}|^p\right)\, .$$
                
Now we introduce the Riemann-Stieljes measure associated to the non-decreasing function $s\mapsto M(s)$ and we have
\begin{equation*}\begin{split}\int_{B(o,R)\setminus B(o,r)} |f|^p|d \chi_{r,R}|^p\,\,&=\int_r^R |\xi'_{r,R}(s)|^p dM(s)\\
&=\varepsilon(r,R)^p \int_r^R  s^{\frac{p}{p-1}} \frac{1}{M(s)^{\frac{p}{p-1}}}dM(s)\\
&\le\varepsilon(r,R)^p\left(r^{\frac{p}{p-1}} \frac{p-1}{M(r)^{\frac{1}{p-1}}}+
\int_r^R  \frac{p\,s^{\frac{1}{p-1}}}{M(s)^{\frac{1}{p-1}}}ds \right)\\
&\le(p-1)\varepsilon(r,R)^p r^{\frac{p}{p-1}} \frac{1}{M(r)^{\frac{1}{p-1}}}+ p\,\varepsilon(r,R)^{p-1}.\end{split}\end{equation*}
With the hypothesis $\int_1^\infty  \left(\frac{r}{M(r)}\right)^{\frac{1}{p-1}}dr=\infty$
it is possible to find two increasing and divergent sequences of $r_\ell<R_\ell$ such that 
$$\lim_{\ell\to \infty} \varepsilon(r_\ell,R_\ell)=0$$ and
$$\lim_{\ell\to \infty} \varepsilon(r_\ell,R_\ell)^p r_\ell^{\frac{p}{p-1}} \frac{1}{M(r_\ell)^{\frac{1}{p-1}}}=0\, .$$
\noindent Hence the result.

\endproof

Assume now that $(M,g)$ is non parabolic and that the Riesz transform is bounded on $L^p$ and on $L^{\frac{p}{p-1}}$. The Hodge projector
$$\Pi=d\Delta^{-1}d^*=(d\Delta^{-\frac12}) (\Delta^{-\frac12}d^*)=(d\Delta^{-\frac12})(d\Delta^{-\frac12})^*$$ extends from 
$L^2\cap L^p$ to a bounded operator on 
$L^p(T^*M)$.  Hence $\Pi(\cC^\infty_0(T^*M))$ is dense in $\Pi(L^p(T^*M))$. 
Also by definition $\Pi$ is the identity on $d\cC^\infty_0(M)$ hence
$$\overline{d\cC^\infty_0(M)}^{L^p}\subset \Pi(L^p(T^*M)).$$
If $\alpha\in \cC^\infty_0(T^*M)$ then we have $\Pi(\alpha)=df$ where
$$f(x)=\int_M G(x,y)d^*\alpha(y) dy$$
and where $G(x,y)$ is the Green kernel of $(M,g)$. Because $d^*\alpha$ has compact support, the growth of $r\mapsto \int_{B(o,r)} f^p$ is controlled by
$$\cG_p(r):=\int_{B(o,r)} G(x,o)^p dx.$$
\noindent Hence a direct consequence of the last lemma is the following proposition:
\begin{prop}\label{prop:cond1}
If $(M,g)$ is a complete Riemannian manifold such that 
\begin{itemize}
 \item $(M,g)$ is non parabolic and its Green kernel satisfies:
$$\int_1^\infty \left(\frac{r}{\cG_p(r)}\right)^{\frac{1}{p-1}}dr=\infty.$$
\item The Riesz transform is bounded on $L^p$ and on $L^{\frac{p}{p-1}}$
\end{itemize}
then $$\Pi(L^p(T^*M))=\overline{d\cC^\infty_0(M)}^{L^p}.$$
\end{prop}

Recall that $(M,g)$ is said to be $p$-parabolic if we can find a sequence of smooth function with compact support $(\chi_k)$ such that 
$$0\le \chi_k\le 1\,\,\lim_{k\to\infty}\|d\chi_k\|_{L^p}=0\,\,\mathrm{and}\,\, \chi_k\to 1 \mathrm{\,\, uniformly\,\, on\,\,  compact \,\, set}.$$
A consequence of the definition is that on a $p$-parabolic manifold, any bounded function with $L^p$ gradient has its gradient in the $L^p$-closure of $d\cC^\infty_0(M)$. 

\begin{rem}\label{Greenborne} If $(M,g)$ is non-parabolic (i.e. $2$-hyperbolic) then its Green kernel $G(x,y)$ is bounded outside its pole: that is to say if $r>0$  then $x\in M\setminus B(y,r)\mapsto G(x,y)$ is positive and bounded by $\max_{x\in \partial B(y,r)} G(x,y)$. 
\end{rem}In particular we get 
\begin{prop}\label{prop:cond2}
 Assume that $(M,g)$ is a complete Riemannian manifold such that 
\begin{itemize}
 \item $(M,g)$ is non $2$-parabolic and $p$-parabolic for some $p>2$,
\item the Riesz transform is bounded on $L^p$ and on $L^{\frac{p}{p-1}}$
\end{itemize}
then $$\Pi(L^p(T^*M))=\overline{d\cC^\infty_0(M)}^{L^p}.$$
\end{prop}
Those two results should be compared with the one of \cite[lemma 7.1]{CCH} where a Sobolev inequality was assumed.  Then we can show the following adaptation of \cite[corollary 7.5]{CCH}
\begin{thm} Under the hypotheses of the \pref{prop:cond2}, $(M,g)$ has only one end.
\end{thm}
\proof
If $M$ has at most two  ends, we can find a compact set $K\subset M$ such that $M\setminus K=\cO_-\cup \cO_+$ with
$\cO_-\cap \cO_+=\emptyset$ and such that both $\cO_\pm$ are unbounded. And we can build a smooth function $\varphi$ 
such that $\varphi=\pm 1$ on $\cO_\pm$. Then $\Delta\varphi\in \cC^\infty_0(M)$, $d\varphi\in \cC^\infty_0(T^*M)$,
We can defined $h\colon M\rightarrow [-1,1]$ by:
$$h(x)=\varphi(x)-\int_{M} G(x,y)\Delta\varphi(y)dy\,\,,$$
$h$ is a harmonic function and with the \rref{Greenborne}, we know that $h$ is bounded and by construction
$$dh=d\varphi-\Pi(d\varphi).$$
As $(M,g)$ is assumed to be $p$-parabolic, we get that:
$$dh\in \overline{d\cC^\infty_0(M)}^{L^p}\, .$$
So that $\Pi(dh)=dh$; but on $L^2$, we have by construction $\Pi(dh)=0$, hence the contradiction.
\endproof
\begin{cor} Let $p>2$. 
On a non-parabolic and $p$-parabolic manifold with at least two ends, the Riesz transform can not be bounded simultaneously on  $L^p$ and on $L^{\frac{p}{p-1}}$.
\end{cor}
\subsection{Sublinear harmonic function and the Riesz transform}
Our next result shows  that the existence of non constant sublinear harmonic function implies some $L^p$  unboundedness of the Riesz transform.
\begin{prop}Let $(M,g)$ be a complete Riemannian manifold whose Ricci curvature satisfies 
$$\ricci_g\ge -\frac{\kappa^2}{r^2(x)}\, g.$$
Assume moreover that $(M,g)$ is doubling and satisfies the (RCE) condition and that the volume of anchored balls satisfy for some $\mu\ge \nu>2$  and positive constant $c$:
$$ 1<r<R\ \Longrightarrow c\left(\frac{R}{r}\right)^\nu\volu B(o,r)\le \volu B(o,R)\le c R^\mu\,\,.$$
Let $\alpha\in [0,1)$. If $(M,g)$ carries a non constant harmonic function $h$ with $\alpha$-growth: 
$$h(x)=\cO(r^\alpha(x))\,\,,$$
then the Riesz transform is not bounded on $L^p$ for any $\displaystyle p>2$ and $\displaystyle  p\ge\frac{\mu}{1-\alpha}$.
\end{prop}
A quick inspection of the proof below shows that we only need that $\nu>\frac{p}{p-1}$.
\proof By contraposition, we assume that the Riesz transform is bounded on $L^p$ with $p(1-\alpha)\ge \mu$ and consider a harmonic function $h$
such that 
$$h(x)=\cO(r^\alpha(x))$$ and we are going show that necessary $h$ is constant.
We remark that our conditions imply a relative Faber-Krahn inequality, hence by \tref{thm:CD} or \cite{CD},  the Riesz transform is bounded on $L^{\frac{p}{p-1}}.$
The quadratic decay of the negative part of the Ricci curvature together with the Cheng-Yau's gradient estimate implies that 
$$dh(x)=\cO\left(r^{\alpha-1}(x)\right)\,\,,$$
\noindent The volume growth condition  $$\volu B(o,r)=\cO(r^\mu)$$  implies 
$$dh\in L^p.$$ Moreover we also have   
$M(r)=\int_{B(o,r)} h^p\le C r^{\mu+p\alpha}\le Cr^p,$
so that 
$$\left(\int_1^\infty \frac{r}{M(r)}\right)^{\frac{1}{p-1}}dr=\infty, $$
and with the \lref{lem:primi}, we get that: 
$$dh\in  \overline{d\cC^\infty_0(M)}^{L^p}.$$
The volume growth assumption also implies that $(M,g)$ is $p$-parabolic and the \pref{prop:cond1} yields:
$$\Pi(dh)=dh.$$
Let $\alpha\in \cC^\infty_0(T^*M)$ then the Hodge projector being bounded on $L^p$ and on $L^{\frac{p}{p-1}}$, we get 
$$\langle dh,\alpha\rangle=\langle \Pi(dh),\alpha\rangle=\langle dh,\Pi(\alpha)\rangle\, .$$
But $$ \Pi(\alpha)=df$$ where $f$ is given by 
$$f(x)=\int_M G(x,y)d^*\alpha(y) dy.$$
There is a $R>0$ so that $f$ is harmonic outside $B(o,R)$ and using the Green kernel estimate, we know that $f$ tends to zero at infinity.
We are going to estimate the decay of $f$. Let $q=p/(p-1)$.
Because $\Pi$ is bounded on $L^q$, we have $df\in L^q$.

Let $x\in M\setminus B(o,2R)$, $f$ is harmonic on the remote ball $B(x,r(x)/2)$ and $u=|df|$ satisfies the elliptic inequality:
$$\Delta u\le \frac{C}{r(x)^2} u\,\, \mathrm{on}\,\,B(x,r(x)/2)\,\,.$$
The lower bound on the Ricci curvature implies that 
$$u^q(x)\le \frac{C}{\volu B(x,r(x)/2)} \int_{B(x,r(x)/2)} u^q=\frac{o(1)}{\volu B(o,r(x))}.$$
hence:
$$|df(x)|\le \frac{o(1)}{(\volu B(o,r(x)))^{\frac1q}}\,\,.$$
Using the (RCE) condition, we can integrate this inequality along a path starting from $x$ and escaping to infinity and get 
$$|f(x)|\le  \frac{o(r(x))}{(\volu B(o,r(x)))^{\frac1q}}+o(1)\int_{r(x)}^\infty \frac{1}{(\volu B(o,s))^{\frac1q}}ds.$$
Using the reverse doubling condition (RD${}_{\nu}$) for the anchored balls we get that 
$$|f(x)|\le C\frac{o(r(x))}{(\volu B(o,r(x)))^{\frac1q}}.$$
Now we define a function $\chi_k$  by 
$$\chi_k(x)=\begin{cases}
1&\,\,\mathrm{on}\,\, B(o,k)\\
2-\frac{r(x)}{k}  &\,\,\mathrm{on}\,\, B(o,2k)\setminus B(o,k)\\
0&\,\,\mathrm{on}\,\, M\setminus B(o,2k).
\end{cases}$$
Because $dh\in L^p$ and $df\in L^q$, we have 
$$\langle dh,\alpha\rangle=\langle dh,df\rangle=\lim_{k\to \infty} \langle dh,\chi_kdf \rangle=\lim_{k\to \infty} \langle dh,d(\chi_kf)\rangle-\langle fdh,d\chi_k\rangle.$$
However $h$ is harmonic hence 
$$ \langle dh,d(\chi_kf)\rangle=0\,\,,$$
moreover 
\begin{equation*}\begin{split}|\langle fdh,d\chi_k\rangle|&\le C \frac{ \volu B(o,2k)\,o(k)}{(\volu B(o,k))^{\frac1q}}k^{(\alpha-1)}k^{-1}\\
&\le (\volu B(o,k))^{\frac{1}{p}}k^{\alpha-1}o(1)\\
&\le  Ck ^{\frac{\mu}{p}}k^{\alpha-1}o(1)\,\,.\end{split}\end{equation*}
Our hypotheses imply that this quantity tends to zero when $k$ tends to infinity. Eventually we obtain that 
for all  $\alpha\in \cC^\infty_0(T^*M)$:
$\langle dh,\alpha\rangle=0$ hence $dh=0$.
\endproof

In fact, the $\alpha$-H\"older Elliptic
 estimates imply that a sublinear harmonic function with $\beta$-growth and $\beta<\alpha$ is necessary constant.
\begin{lem}\label{Sublne} Assume that  $(M,g)$ carries a non constant sublinear harmonic function $h$ with
$$h(x)=\cO\left(r^\beta(x)\right)\ .$$ If $\alpha>\beta$
then $(M,g)$ can not satisfy the $\alpha$-H\"older Elliptic
 estimates . 
 \end{lem}
\proof Indeed, assume that $(M,g)$ satisfies  (EH${}_{\alpha}$) and consider a harmonic function $h\colon M\rightarrow \R$   such that for some positive constants $C$ and $\beta<\alpha$:
 $$\forall x\in M,\ |h(x)|\le C\left(1+d(x,o)^\beta\right).$$
 Using (EH${}_{\alpha}$) we get that for any $x\in B(o,R)$:
 $$\left|h(x)-h(o)\right|\le \Lambda \left(\frac{d(x,o)}{R}\right)^\alpha\sup_{y\in B(o,2R)} |h(y)|.$$
 Hence we get that for any $x\in M$ and any $R\ge r(x)$:
  $$\left|h(x)-h(o)\right|\le C\Lambda \left(\frac{d(x,o)}{R}\right)^\alpha\left(1+(2R)^\beta\right).$$
 Letting $R\to +\infty$, we get that $h(x)=h(o)$ for any $x\in M$.
 \endproof

\section{Examples}In this section, we describe three series of applications of our results. \subsection{Manifolds with conical ends}\label{explconic} These manifolds $(M^n,g)$ are isometric outside a compact set to a truncated cone:
$$\cC_R(\Sigma):=\,\left((R,\infty)\times\Sigma, (dr)^2+r^2h\right)$$ where 
$(\Sigma,h)$ is compact Riemannian manifold. We are going to explain how our results can be used to recover the \tref{Conic}  of H-Q. Li, C. Guillarmou and A. Hassell  and of A. Hassell and P. Lin (\cite{Li,GH1, GH2,HP}).
From the explicit form of the metric, it is easy to remark that the conditions (QD), (VC) and  (RCE) are satisfied. Concerning the volume growth of geodesic balls, there is a positive constant $C$ such that for all $x\in M$ and $R>0$:
$$C^{-1} R^n\le \volu B(x,R)\le C R^n\ .$$ 
\begin{itemize}
\item[$\star$]{\bf A general positive result:} The hypotheses of \tref{thm:general} are satisfied for $\nu=n$, hence
the Riesz transform
is bounded on $L^p$ for all $p\in (1,n)\cup\{2\}$.

\item[$\star$] {\bf Negative results if $M^n$ has several ends (i.e. if $\Sigma$ is not connected):} In the case where $n>2$, \cite[Corollary 7.5]{CCH} already told us that the Riesz transform is  unbounded on $L^p$ when $p\ge n$. In the case where $n=2$, it is well known that $(M,g)$ carry a non constant harmonic function with logarithmic growth. Indeed $(M^2,g)$ is conformally equivalent to a compact surface  with a finite number of points removed:
$$(M^2,g)=\left(\Sigma\setminus \{p_1,\dots,p_r\},\varphi^{-2}\bar g\right)$$ where $(\Sigma,\bar g)$ is  compact Riemann surface. Moreover around each $p_i$ we have 
$$d_{\bar g}(x,p_i)\simeq d_g(x,o)^{-1}.$$
Consider $\bar G$ the Green kernel on $(\Sigma,\bar g)$, then the function 
$f(x)=G(x,p_1)-G(x,p_2)$ is harmonic on $\left(\Sigma\setminus \{p_1,p_2\},\bar g\right)$ and satisfies for $i=1,2$:
$$f(x)\simeq (-1)^{i} \log\left((d_{\bar g} (x,p_i)\right).$$  Recall $\Delta_g=\varphi^{2}\Delta_{\bar g}$, hence $f$ is a harmonic function on $(M^2,g)$ and it has logarithmic growth:
$f(x)=\cO\left(\log r(x)\right)$. The \pref{subNO} imply that if $p>2$ then the Riesz transform can not be bounded in $L^p$ .

 Assume now that $(\Sigma,h)$ is connected.  Let $$0<\lambda_1\le\lambda_2\le \dots$$ be the spectrum of the Laplace operator of $(\Sigma,h)$, where the eigenvalues are repeated according to their multiplicity. 
 And let $\varphi_0,\varphi_1\dots$ be an associated set of normalized eigenfunctions:
 $\Delta_h \varphi_i=\lambda_i\varphi_i$. We define
$0<\alpha_1\le \alpha_2\dots $ by 
$$\alpha_i(n-2-\alpha_i)=\lambda_i\ .$$
\item[$\star$] {\bf Negative results if $M^n$ has only one end:} According to \cite[Lemma 6]{CZw}, $(M,g)$ carries a harmonic function $f$ such that 
$$f(r,\theta)=r^{\alpha_1}\varphi_1(\theta)+o\left(r^{\alpha_1}\right).$$ 
Define now
$$\alpha=\min(\alpha_1,1).$$
The  \pref{subNO} imply that if $p(1-\alpha)\ge n$, then the Riesz transform can not be bounded in $L^p$.
\item[$\star$] {\bf Positive results if $M^n$ has only one end:} Assume that $(\Sigma,h)$ is connected.
Let $\cB(r)=K\cup\left([1,r]\times \Sigma\right)$. 
Following the analysis of \cite[Proposition 4.1]{ACM}, it can be shown that there is a constant $C$ such that for all harmonic function $f$ defined over $\cB(R)$ and all 
$1<r<\rho\le R$:
$$\frac{r^{2-2\alpha}}{r^n}\int_{\cB(r)} |df|^2\le C\frac{\rho^{2-2\alpha}}{\rho^n}\int_{\cB(\rho)} |df|^2.$$
From this result and \pref{prop:freq}, we obtain that $(M,g)$ satisfies the $\alpha$-H\"older estimates (EH${}_{\alpha}$). Hence if $p(1-\alpha)<n$, then the Riesz transform is bounded on $L^p$ .
\end{itemize}

\subsection{Model manifolds}\label{WPsection} We consider a Riemannian manifold $(M^n,g)$  that is isometric outside a compact set to the warped product:
$$\cC_R(\Sigma):=\,\left((R,\infty)\times\Sigma, (dr)^2+f^2(r)h\right)$$ where 
$f(r)=e^{u(\ln(r))},$
for some function $u\colon [1,\infty)\rightarrow \R$ with bounded second derivative and
$$\alpha\le u'\le 1. $$
We assume that $(\Sigma,h)$ is a compact manifold with non negative Ricci curvature so that $(M,g)$ satisfies the quadratic decay  on the negative part 
of the Ricci curvature and the conditions (VC) \& (RCE). Moreover anchored balls satisfies the reverse doubling condition (RD${}_{\nu}$) for the exponent $\nu=(n-1)\alpha+1$. \par 
If $u$ satisfies the asymptotic condition:
\begin{equation}
 \label{integrable} u'(\theta)\ge \bar\alpha-\psi(\theta) \mathrm{\,\, for\,\, some\,\,non\,\, negative\ }\psi\in L^1
\end{equation}
then the exponent in the reverse doubling is improved to $\bar\nu=(n-1)\bar\alpha+1$. Hence the \tref{thm:general} implies  that the Riesz transform is bounded on $L^p$ for $p<(n-1)\bar\alpha+1$.

If $\int_1^\infty \frac{1}{f^{n-1}(r)}dr<\infty\,$ then $(M,g)$ is non-parabolic and if $\int_1^\infty \frac{1}{f^{\frac{n-1}{p-1}}(r)}dr=\infty\,$ then $(M,g)$ is $p$-parabolic. Hence by \tref{thm:nob}: if $p>2$ and if  $\Sigma$ is non connected then the Riesz transform is not bounded on $L^p$.

If $\Sigma$ is connected and if $\int_1^\infty \left(1-u'(r)\right)dr=\infty$ then the diameter of geodesic sphere growth slowly
$$\diam \partial B(o,R)= \sup_{x,y\in \partial B(o,R)} d(x,y)=o(R)$$ and we can apply the \cref{noharm} and we know that the Riesz transform is bounded on $L^p$ for all $p\in (1,+\infty)$.
In conclusion, we get:
 \begin{prop} \begin{itemize}
 \item If the function $u$ satisfies the condition (\ref{integrable}) and if $p<(n-1)\bar\alpha+1$  then the Riesz transform is bounded on $L^p$. 
 \item Assume that $1<(n-1)\bar\alpha$ and that $\Sigma$ is non connected. If
$$\int_1^\infty \frac{1}{f^{\frac{n-1}{p-1}}(r)}dr=\infty\,$$
then the Riesz transform not bounded on $L^p$. 
\item If $\Sigma$ is connected and if $\int_1^\infty \left(1-u'(r)\right)dr=\infty$ then the Riesz transform is bounded on $L^p$ for all $p\in (1,+\infty)$.

\end{itemize}
\end{prop}

\subsection{Examples with infinite topological type} Our last example is inspired by a construction of J.Lott and Z.Shen \cite{LottShen}. 
Assume that $(P^n,g)$ is a compact Riemannian manifold with boundary $\partial P=\Sigma_-\cup \Sigma_-$ and that $\Sigma_\pm$ have collared neighborhoods
$\cU_\pm$ and assume moreover that there is a diffeomorphism
$f\colon \cU_-\rightarrow \cU_+$ with $$f^*g=4g.$$
For $\ell\in \N$, we define $2^\ell P$ to be the rescaled Riemannian manifold
$$(P,4^\ell g)\,.$$ 
Using the map $f$ we can glue all the  $2^\ell P$ and get a Riemannian manifold $(X,g)$ with boundary 
$\partial X=\Sigma_-$.
 If there is a compact manifold $K$ with boundary diffeomorphic to $\Sigma_-$ then we can form
$X_0=X\cup K$. We can also form $X_1=X\#_{\Sigma_-}(-X)$ the double of $X$.
These manifolds have quadratic curvature decay and Euclidean growth. Hence our results yield the following:
\begin{prop}\begin{itemize}
             \item On $X_0$ and $X_1$, the Riesz transform is bounded on $L^p$ for every $p\in (1,n)$.
\item On $X_1$, the Riesz transform is not bounded on $L^p$ if $p\ge n>2$.
\item Assume that $P$ is connected, then there is some $\varepsilon\in (0,1]$ such that on $X_0$, the Riesz transform is $L^p$-bounded for every $p\in 
\left(1,\frac{n}{1-\varepsilon}\right)$.
            \end{itemize}

\end{prop}

\section{Some perpectives}
In this last section, we conclude by some remarks and perpectives:
\subsection{Perturbations}
In \cite{Coulhon-Dungey}, T. Coulhon and N. Dungey obtained a nice property of stability under perturbation for $L^p$ boundedness of the  Riesz transform. This result implies that  if $(M,g)$ is a complete Riemannian manifold satisfying the hypotheses of the \tref{thm:general} (resp. \tref{thm:impro}) and a non collapsing hypothesis:
$$\inf_{x\in M} \volu B(x,1)\,>0$$ then
the same conclusions hold for any other metric $\tilde g$ that satisfies for some $\upvarepsilon>0$: 
$$\tilde g-g=\cO\left(r^{-\upvarepsilon}(x)\right).$$

\subsection{Riesz transform associated to Schr\"odinger operator}
Let $(M,g)$ be a complete Riemannian manifold satisfying the geometric condition (QD), (VC) and (RCE) and $V\colon M\rightarrow \R$ be a locally bounded function that satisfies 
$$V(x)=\cO\left(r^{-2}(x)\right).$$
A challenging question is to study the boundedness of the Riesz transform associated to Schr\"odinger operator $\Delta+V$:
$$ R_V=d(\Delta+V)^{-\frac12}\, .$$ 
In this case, the situation can be very complicated even in the case of manifolds with conical ends (see \cite{GH1,GH2, HP}). Some part of our analysis extends easily to this case but some crucial points are missing, we leave as an open question of finding the appropriate spectral condition for the operator $\Delta+V$ that leads to results similar to  \tref{thm:general} or \tref{thm:impro}. Note that in the case where the potential is small in some Kato class, then the $L^p$-boundedness of the Riesz transform $R_0$ implies the  $L^p$-boundedness of the Riesz transform $R_V$ (\cite{DevyverP}) .
\subsection{Gaussian estimates on the heat kernel on $1$-forms}
In the setting of \tref{thm:impro}, it will be interesting to understand under which geometric assumptions, one get a Gaussian estimate on the heat kernel of the Hodge Laplacian on $1$-forms. Recall that such an estimate implies a boundedness result for the Riesz transform (\cite{CD2,CDfull,siko}). According to \cite{CZ,DA,DevyverP}, a sub-critical assumption on the Ricci curvature implies such an upper bound. It has recently been shown that this sub-critical assumption yields results for the Riesz transform on $1$ and $2$-forms (\cite{Magniez}).

\subsection{Riesz transform associated to differential forms}
Let $\Delta_k=dd^*+d^*d\colon \cC^\infty_0(\Lambda^kT^*M)\rightarrow \cC^\infty_0(\Lambda^kT^*M)$ be the Hodge Laplacian acting on $k$-differential forms. Following the questions asked in \cite{CCH}, one would like to understand the $L^p$ boundedness of the Riesz transforms
$$d(\Delta_k)^{-\frac12},\, d^*(\Delta_k)^{-\frac12}, \nabla(\Delta_k)^{-\frac12}.$$
On manifold with conical ends, these questions have been recently investigated by C. Guillarmou and D. Sher  (\cite{GSher}) and it  appears that it is a very difficult problem; it is tempting to analyze what can be done on a manifold whose curvature tensor decays quadratically, however in order to use our analysis, we will need Gaussian upper bound on the heat kernel on forms and such an estimate already provided some boundedness result for the Riesz transform according to a general principle (\cite{CD2,siko}).
\subsection{The Riesz transform of second order.}
Another interesting question is about the $L^p$ boundedness of the Riesz transform of second order:
$$\nabla d\Delta^{-1}\,.$$
The study of the $L^p$ boundedness of the operator 
$\nabla d(\Delta+1)^{-1}\,$ has been recently investigated by  B. G\"uneysu and S. Pigola (\cite{BS}). 
Motivated by this paper, 
in a future work, we intend to consider applications of these ideas for the operator $\nabla d\Delta^{-1}$ on manifold where the full curvature tensor decays quadratically:
$$\|Rm(x)\|=\cO\left(\frac{1}{r^2(x)}\right).$$

\end{document}